\documentclass[reqno,12pt,a4paper]{amsart}
\usepackage{amscd,amsfonts,amssymb,amsthm,latexsym}
\usepackage{srcltx}

\theoremstyle{plain}
\newtheorem{theorem}{Theorem}[section]
\newtheorem*{theorem*}{Theorem}
\newtheorem{corollary}{Corollary}[section]
\newtheorem{lemma}{Lemma}[section]

\theoremstyle{definition}

\newtheorem*{definition*}{Definition}

\theoremstyle{remark}

\newtheorem*{remark*}{Remark}

\numberwithin{equation}{section}

\begin{document}
\raggedbottom 

\title[On Romanoff's theorem]{On Romanoff's theorem}

\author{Artyom Radomskii}

\begin{abstract} We obtain some results related to Romanoff's theorem.
\end{abstract}

 \address{HSE University, RF}

\thanks{The article was prepared within the framework of the Basic Research Program at HSE University, RF}

\email{artyom.radomskii@mail.ru}

\maketitle

\section{Introduction}
Let $\varphi$ denote Euler's totient function. It is clear that $1\leq \varphi(n)\leq n$ for any positive integer $n$. Therefore, if $a_1,\ldots, a_{N}$ are positive integers (not necessarily distinct), $s\in \mathbb{N}$, then
 \[
\sum_{n=1}^{N} \biggl(\frac{a_{n}}{\varphi(a_{n})}\biggr)^{s}\geq N.
\] Upper bounds for such sums are of interest. It is well-known, for any positive integer $s$ there is a positive constant $c(s)$, depending only on $s$, such that
 \[
 \sum_{n=1}^{N} \left(\frac{n}{\varphi (n)}\right)^{s}\leq c(s) N
 \] for any positive integer $N$.

  We prove the following result.
\begin{theorem}\label{T1}
Let $\alpha$ be a real number with $0<\alpha <1$. Then there is a constant $C(\alpha)>0$, depending only on $\alpha$, such that the following holds. Let $M$ be a real number with $M\geq 1$, let $a_1,\ldots, a_{N}$ be positive integers \textup{(}not necessarily distinct\textup{)} with $a_{n}\leq M$ for all $1\leq n \leq N$. We define
\[
\omega(d)=\#\{1\leq n\leq N: a_{n}\equiv 0\ \textup{\text{(mod $d$)}}\}
\]for any positive integer $d$. Let $s$ be a positive integer. Then
\[
\sum_{n=1}^{N}\left(\frac{a_{n}}{\varphi(a_{n})}\right)^{s}\leq (C(\alpha))^{s}
\biggl(N+\sum_{p\leq (\ln M)^{\alpha}} \frac{\omega(p)(\ln p)^{s}}{p}\biggr).
\]
\end{theorem}

 One can see from the proof of Theorem \ref{T1} that $C(\alpha)=c/\alpha$, where $c$ is a positive absolute constant. The following result shows that Theorem \ref{T1} can not be improved in the following direction: a condition $p\leq (\ln M)^{\alpha}$ can not be replaced by a condition $p\leq (\ln M)^{o(1)}$.

\begin{theorem}\label{T2}
Let $\alpha(M)$, $M=1, 2, \ldots,$ be a sequence of positive real numbers such that $\alpha(M)\to 0$ as $M\to +\infty$ and $(\ln M)^{\alpha(M)}\geq 2$ for all $M\geq 3$. Then there is a constant $M_0 >0$, depending only on a sequence $\alpha(M)$, such that the following holds. For any positive integer $M\geq M_0$ there is a non-empty set $A\subset \{1,\ldots, M\}$ such that
\[
\#\{n\in A: n\equiv 0\ \text{\textup{(mod $p$)}}\}=0
\]for any prime $p\leq (\ln M)^{\alpha(M)}$ and
\[
\sum_{n\in A}\frac{n}{\varphi(n)}\geq \frac{c}{\alpha(M)}\#A.
\]Here $c>0$ is an absolute constant.
\end{theorem}

From Theorem \ref{T1} we obtain
\begin{theorem}\label{T3}
Let $\varepsilon$ be a real number with $0<\varepsilon <1$. Then there is a constant $C(\varepsilon)>0$, depending only on $\varepsilon$, such that the following holds. Let $x$ and $z$ be real numbers with $x\geq 3$ and $(\ln x)^{\varepsilon} \leq z \leq x$. Let $k$ be a positive integer, let $a_0,\ldots, a_k$ be integers with $|a_i| \leq x$ for all $0\leq i \leq k$, and $a_k \neq 0$. By $\delta:= (a_0,\ldots, a_k)$ we denote the greatest common divisor of $a_0,\ldots, a_k$. Let
\[
R(n)=a_k n^{k}+a_{k-1}n^{k-1}+\ldots+a_0.
\]Let $s$ be a positive integer. Then
\begin{equation}\label{Ineq.Th2}
\sum_{\substack{-z\leq n \leq z\\ R(n)\neq 0}}\left(\frac{|R(n)|}{\varphi(|R(n)|)}\right)^{s}
\leq \biggl(C(\varepsilon)\frac{\delta}{\varphi(\delta)} \ln (k+1) \biggr)^{s}s!\, z.
\end{equation}
\end{theorem}

\begin{corollary}\label{Cor1}
 Let $k$ be a positive integer, let
\[
R(n)=a_k n^k+a_{k-1}n^{k-1}+\ldots+a_0
\] be a polynomial with integer coefficients, $a_k\neq 0$. Then there is a constant $C(R)>0$, depending only on $R$, such that if $s$ is a positive integer and $x$ is a real number with $x\geq 1$, then
\[
\sum_{\substack{-x\leq n \leq x\\ R(n)\neq 0}}\left(\frac{|R(n)|}{\varphi(|R(n)|)}\right)^{s}
\leq (C(R))^{s}s!\, x.
\]
\end{corollary}

Let $\mathcal{L}=\{L_1,\ldots,L_k\}$ be a set of $k$ linear functions with integer coefficients
\[
L_{i}(n)=a_i n+b_i,\qquad i=1,\ldots, k.
\]For $L(n)=an+b$, $a, b\in \mathbb{Z}$, we define
\[
\Delta_{L}=|a|\prod_{i=1}^{k}|a b_i - b a_i|.
\]In modern sieve methods sums
\[
\sum_{(a,b)\in \Omega} \frac{\Delta_{L}}{\varphi(\Delta_{L})}
\]appear (see, for example, \cite{Maynard}). Here $(a,b)$ denotes a vector and $\Omega$ is a finite set in $\mathbb{Z}^2$. From Theorem \ref{T1} we obtain
\begin{theorem}\label{T4}
 Let $\varepsilon$ be a real number with $0<\varepsilon <1$. Then there is a constant $C(\varepsilon)>0$, depending only on $\varepsilon$, such that the following holds. Let $x$ and $z$ be real numbers with $x\geq 3$ and $(\ln x)^{\varepsilon} \leq z \leq x$. Let $a, b_1,\ldots, b_k$ be integers with $a\geq 1$, $|b_i| \leq x$ for all $1\leq i \leq k$. Let $\mathcal{L}=\{L_1,\ldots,L_k\}$ be a set of $k$ linear functions, where $L_{i}(n)=a n+b_i,\ i=1,\ldots, k$. For $L(n)=an+b$, $b\in \mathbb{Z}$, we define $\Delta_{L}=a^{k+1}\prod_{i=1}^{k}|b_i - b|$. Let $s$ be a positive integer. Then
\begin{equation}\label{Ineq.Th3}
\sum_{\substack{-z\leq b \leq z\\ L(n)=an+b \notin \mathcal{L}}}\left(\frac{\Delta_{L}}{\varphi(\Delta_L)}\right)^{s}
\leq \biggl(C(\varepsilon)\frac{a}{\varphi(a)} \ln (k+1) \biggr)^{s}s!\, z.
\end{equation}
\end{theorem}

Theorem \ref{T4} extends a result of Maynard \cite[Lemma 8.1]{Maynard} which showed the same result but with $s=1$ and $x^{1/10} \leq z \leq x$. Since $a/\varphi(a) \leq c\ln\ln (a+2)$, where $c$ is a positive absolute constant, the right-hand side of \eqref{Ineq.Th3} can be replaced by
\[
\bigl(C(\varepsilon)\ln\ln (a+2) \ln (k+1) \bigr)^{s}s!\, z.
\] The analogous remark is true for \eqref{Ineq.Th2}.

We recall some facts on elliptic curves (for more details, see, for example, \cite[Chapter XXV]{Hardy.Wright}). An \emph{elliptic curve} is given by an equation of the form
\[
E: y^2 = x^3 + Ax+B,
\]with the one further requirement that the \emph{discriminant}
\[
\Delta= 4A^3+ 27 B^2
\]should not vanish. For convenience, we shall generally assume that the coefficients $A$ and $B$ are integers. One of the properties that make an elliptic curve $E$ such a fascinating object is the existence of a composition law that allows us to 'add' points to one another. We adjoin an idealized point $\mathcal{O}$ to the plane. This point $\mathcal{O}$ is called the \emph{point at infinity}. The special rules relating to the point $\mathcal{O}$ are
\[
P+(-P)=\mathcal{O}\quad \text{and}\quad P+\mathcal{O}=\mathcal{O}+P=P
\]for all points $P$ on $E$. Given a prime $p$, by $\mathbb{F}_{p}$ we denote the field of classes of residues modulo $p$. We define
\[
E(\mathbb{F}_{p})=\{(x,y)\in {\mathbb{F}}_{p}^2: y^2\equiv x^3+Ax+B\ \text{(mod $p$)}\}\cup\{\mathcal{O}\}.
\]Repeated addition and negation allows us to 'multiply' points of $E$ by an arbitrary integer $m$. This function from $E$ to itself is called the \emph{multiplication-by-$m$} map,
\[
\phi_{m}: E\to E,\qquad \phi_{m}(P)=mP=\text{\textup{sign}}(m)(P+\dots+P)
\](the sum contains $|m|$ terms). By convention, we also define $\phi_{0}(P)=\mathcal{O}$. The multiplication-by-$m$ map is defined by rational functions. Maps $E\to E$ defined by rational functions and sending $\mathcal{O}$ to $\mathcal{O}$  are called \emph{endomorphisms} of $E$. For most elliptic curves (over the field of complex numbers $\mathbb{C}$), the only endomorphisms are the multiplication-by-$m$ maps. Curves that admit additional endomorphisms are said to have \emph{complex multiplication}.

Let $\pi(x)$ denote the number of primes not exceeding $x$. We prove
\begin{theorem}\label{T5}
Let $E$ be an elliptic curve given by an equation
\[
y^2=x^3+Ax+B,
\] where $A$ and $B$ are integers with $\Delta=4A^3+27B^2\neq 0$. Suppose that $E$ does not have complex multiplication. Let $s$ be a positive integer and $x$ be a real number with $x\geq 2$. Then
\begin{equation}\label{T5.Ellip.CURVE.INEQ}
\pi(x)\leq \sum_{p\leq x}\biggl(\frac{\#E(\mathbb{F}_{p})}{\varphi(\#E(\mathbb{F}_{p}))}\biggr)^{s}
\leq C(E,s)\pi(x),
\end{equation}where $C(E,s)>0$ is a constant, depending only on $E$ and $s$.
\end{theorem}

Let $\mathbb{P}$ denote the set of all prime numbers. In 1934 Romanoff proved the following result.\\[-4pt]

\textsc{Theorem (Romanoff \cite{Romanoff}).} \emph{Let $a$ be an integer with  $a\geq 2$. Then there is a constant $c(a)>0$, depending only on $a$, such that
\[
\#\{1\leq n \leq x: \text{there are $p\in \mathbb{P}$ and $j\in \mathbb{Z}_{\geq 0}$}\text{ such that $p+ a^{j}=n$}\}\geq c(a)x
\]for any real number $x\geq 3$.}\\[-5pt]

We prove

\begin{theorem}\label{T6}
Let $A=\{a_{n}\}_{n=1}^{\infty}$ be a sequence of positive integers \textup{(}not necessarily distinct\textup{)}. We put
\begin{align*}
N_{A}(x)&=\#\{j\in \mathbb{N}: a_j \leq x\},\\
\textup{ord}_{A}(n)&=\#\{j\in \mathbb{N}: a_j=n\},\qquad  n\in\mathbb{N},\\
\rho_{A}(x)&=\max_{n\leq x}\textup{ord}_{A}(n).
\end{align*}Suppose that $\textup{ord}_{A}(n)<+\infty$ for any positive integer $n$. Suppose that there are constants $\gamma_{1}>0$, $\gamma_{2}>0$, $\alpha>0$, $x_{0}\geq 10$ such that the following holds.  For any real number $x\geq x_0$, we have
\begin{gather}
N_{A}(x)>0,\label{T6.INEQ1}\\ N_{A}\biggl(\frac{x}{2}\biggr)\geq \gamma_1 N_{A}(x),\label{T6.INEQ2}\\
\sum_{\substack{k\in \mathbb{N}:\\ a_k < x}}\sum_{p\leq (\ln x)^{\alpha}}
\frac{\#\{j\in\mathbb{N}: a_{k}< a_j \leq x\ \text{and } a_j\equiv a_k\ \textup{(mod $p$)}\}\ln p}{p}\leq \gamma_2 (N_{A}(x))^2.\label{T6.INEQ3}
\end{gather}For any positive integer $n$, we put
\[
r(n)=\#\{(p,j)\in \mathbb{P} \times\mathbb{N}: p+a_j=n\}.
\]Then there are constants $c_1 = c_{1}(\gamma_{1})>0$ and $c_2=c_{2}(\gamma_1, \gamma_2, \alpha)>0$, depending only on $\gamma_1$ and $\gamma_1$, $\gamma_2$, $\alpha$ respectively, such that
\[
\#\biggl\{1\leq n \leq x: r(n)\geq c_{1}\frac{N_{A}(x)}{\ln x}\biggr\}\geq c_{2}x\frac{N_{A}(x)}{N_{A}(x)+\rho_{A}(x)
\ln x}
\]for any real number $x\geq x_0$. In particular,
\begin{align*}
\#\{1&\leq n \leq x: \text{there are $p\in \mathbb{P}$ and $j\in \mathbb{N}$ such that $p+a_j = n$}\}\\
&\geq c_{2}x\frac{N_{A}(x)}{N_{A}(x)+\rho_{A}(x)\ln x}
\end{align*}for any real number $x\geq x_0$.
\end{theorem}
We note that Romanoff's theorem follows from Theorem \ref{T6}.

 From Theorem \ref{T6} we obtain
 \begin{theorem}\label{T7}
Let $k$ be an integer with $k\geq 2$ and
\[
R(n)=a_k n^k+\ldots+a_0
\]be a polynomial with integer coefficients, $a_k>0$. For any positive integer $n$, we put
\[
r(n)=\#\{(p, j)\in \mathbb{P}\times\mathbb{N}: p+R(j)=n\}.
\]Then there are constants $c_1>0$, $c_2>0$, $x_0>0$, depending only on $R$, such that
\[
\#\biggl\{1\leq n \leq x: r(n)\geq c_1 \frac{x^{1/k}}{\ln x}\biggr\}\geq c_2 x
\]for any real number $x\geq x_0$.
\end{theorem}

\begin{corollary}\label{Cor2}
Let $k$ be an integer with $k\geq 2$. For any positive integer $n$, we put
\[
r(n)=\#\{(p,j)\in \mathbb{P}\times\mathbb{N}: p+j^k=n\}.
\]Then there are constants $c_1(k)>0$ and $c_2(k)>0$, depending only on $k$, such that
\[
\#\biggl\{1\leq n \leq x: r(n)\geq c_1(k) \frac{x^{1/k}}{\ln x}\biggr\}\geq c_2(k) x
\] for any real number $x\geq 3$. In particular,
\begin{equation}
\#\{1\leq n \leq x: \text{there are $p\in \mathbb{P}$ and $j\in \mathbb{N}$ such that $p+j^k = n$}\}\geq c_2(k)x\label{C2.Roman}
\end{equation}for any real number $x\geq 3$.
\end{corollary}
Corollary \ref{Cor2} extends a result of Romanoff which showed the inequality \eqref{C2.Roman}.

\begin{theorem}\label{T8}
 Let $E$ be an elliptic curve given by an equation $y^2=x^3+Ax+B$, where $A$ and $B$ are integers with $\Delta=4A^3+27B^2\neq 0$. Suppose that $E$ does not have complex multiplication. For any positive integer $n$, we put
\[
r(n)=\#\{(p, q)\in \mathbb{P}^{2}: p+\#E(\mathbb{F}_q)=n\}.
\]Then there are constants $x_0>0$, $c_1>0$, $c_2(E)>0$, where $x_0$ and $c_1$ are absolute constants, $c_2(E)$ is a constant depending only on $E$, such that
\[
\#\biggl\{1\leq n\leq x: r(n)\geq c_1\frac{x}{(\ln x)^{2}}\biggr\}\geq c_{2}(E)x
\]for any real number $x\geq x_0$.
\end{theorem}

\begin{theorem}\label{T9}
 Let $a$ and $b$ be integers with $a\geq 2$ and $b\geq 2$. Then there are positive constants $c_1 (a,b)$ and $c_2(a,b)$, depending only on $a$ and $b$, such that
\begin{align*}
c_1(a&, b)\frac{x}{(\ln x)^{1-1/b}}\\
&\leq \#\{1\leq n \leq x: \text{there are $p\in \mathbb{P}$ and $j\in \mathbb{Z}_{\geq 0}$ such that $p+ a^{j^{b}}=n$}\}\\
 &\leq c_2(a,b)\frac{x}{(\ln x)^{1-1/b}}
\end{align*}for any real number $x\geq 3$.
\end{theorem}

\section{Notation}
We reserve the letters $p$, $q$ for primes. In particular, the sum $\sum_{p\leq K}$ should be interpreted as being over all prime numbers not exceeding $K$. By $\pi(x)$ we denote the number of primes not exceeding $x$. Let $\#A$ denote the number of elements of a finite set $A$. By $\mathbb{Z}$, $\mathbb{Z}_{\geq 0}$, $\mathbb{N}$, $\mathbb{Q}$, $\mathbb{R}$, and $\mathbb{C}$  we denote the sets of all integers, non-negative integers, positive integers, rational numbers, real numbers, and  complex numbers respectively. By $\mathbb{P}$ we denote the set of all prime numbers. Let $(a,b)$ be the greatest common divisor of integers $a$ and $b$, and $[a,b]$ be the least common multiple of integers $a$ and $b$. If $d$ is a divisor of $b-a$, we say that $b$ is congruent to $a$ modulo $d$, and write $b \equiv a$ (mod $d$). Let $\varphi$ denote Euler's totient function, i.\,e. $\varphi(n)=\#\{1\leq m \leq n:\ (m,n)=1\}$. We write $\nu(n)$ for the number of distinct prime divisors of $n$, and $\tau(n)$ for the number of positive divisors of $n$. Let $P^{+}(n)$ denote the greatest prime factor of $n$, and $P^{-}(n)$ denote the least prime factor of $n$ (by convention $P^{+}(1)=1$, $P^{-}(1)=+\infty$). We denote by $\mathcal{M}$ the set of square-free numbers, that is, the number $1$ and the positive integers of the form $p_1\cdots p_l$, where $p_1, \ldots, p_{l}$ are distinct primes. By definition, we put
\[
\sum_{\varnothing} = 0,\qquad \prod_{\varnothing}=1.
\]The symbol $b|a$ means that $b$ divides $a$. For fixed $a$ the sum $\sum_{b|a}$ and the product $\prod_{b|a}$ should be interpreted as being over all positive divisors of $a$. If $x$ is a real number, then $[x]$ denotes its integral part, and $\lceil x\rceil$ is the smallest integer $n$ such that $n\geq x$. We put $\log_{a}x:=\ln x/\ln a$.

For real numbers $x,$ $y$ we also use $(x,y)$ to denote the open interval, and $[x,y]$ to denote the closed interval. Also by $(a_1,\ldots, a_n)$ we denote a vector. The usage of the notation will be clear from the context.

\section{Proofs}
\textsc{Proof of Theorem \ref{T1}.} First we prove the following
\begin{lemma}\label{L1}
Let $n$ be an integer with $n>1$, and $y$ be a positive real number. Then
\[
\prod_{p|n:\, p>y}\biggl(1+\frac{1}{p}\biggr)\leq \textup{exp}\,\frac{\nu(n)}{y}.
\]
\end{lemma}

\textsc{Proof of Lemma \ref{L1}.} We put $\Omega = \{p: p|n\text{ and $p>y$}\}$. Let us consider two cases.

1) Let $\Omega=\varnothing$. Since $n>1$, we have $\nu(n)\geq 1$. We obtain
\[
\prod_{p|n:\,\, p>y}\biggl(1+\frac{1}{p}\biggr)=\prod_{\varnothing} = 1 \leq \textup{exp}\,\frac{\nu(n)}{y}.
\]

2) Suppose that $\Omega\neq\varnothing$. Using the inequality $1+x\leq e^{x}$, $x\in \mathbb{R}$, we have
\[
\prod_{p|n:\,\, p>y}\biggl(1+\frac{1}{p}\biggr)\leq \textup{exp}\biggl(\sum_{p|n:\,\, p>y}\frac{1}{p}\biggr)\leq
\textup{exp}\,\frac{\nu(n)}{y}.
\]Lemma \ref{L1} is proved.

  \begin{lemma}\label{L2}
 Let $n$ be a positive integer. Then
 \[
 \prod_{p|n:\,\, p>\ln n}\biggl(1+\frac{1}{p}\biggr)\leq 5.
 \]
  \end{lemma}

 \textsc{Proof of Lemma \ref{L2}.} If $n=1$, then the product equals 1 and the statement is true. Let $n>1$. It is clear that
 \begin{equation}\label{v(n)<log_n}
 \nu(n)\leq \log_{2}n=\frac{\ln n}{\ln 2}.
 \end{equation}Applying Lemma \ref{L1} with $y=\ln n$ and  \eqref{v(n)<log_n}, we obtain
 \[
\prod_{p|n:\,\, p>\ln n}\biggl(1+\frac{1}{p}\biggr)\leq \textup{exp}\,\frac{\nu(n)}{\ln n}\leq \textup{exp}\,\frac{1}{\ln 2}< 5.
\]Lemma \ref{L2} is proved.

\begin{lemma}\label{L3}
 Let $\alpha$ be a real number with $0<\alpha < 1$. Then there is a constant $C(\alpha)>0$, depending only on $\alpha$, such that if $n$ is a positive integer, then
\[
\frac{n}{\varphi(n)}\leq C(\alpha)\prod_{p|n:\,\,
p\leq (\ln n)^{\alpha}} \biggl(1+\frac{1}{p}\biggr).
\]
\end{lemma}

\textsc{Proof of Lemma \ref{L3}.} We may assume that $n\geq \exp(2^{1/\alpha})$. By $\zeta (s)$ we denote the Riemann zeta function. We have
\begin{align*}
\frac{n}{\varphi(n)}&= \prod_{p|n} \biggl(1-\frac{1}{p}\biggr)^{-1}=\prod_{p|n} \frac{p}{p-1}
=\prod_{p|n} \frac{p}{p-1} \frac{p}{p+1}\frac{p+1}{p}\\
&\leq \prod_{p|n} \biggl(1+\frac{1}{p}\biggr)\prod_{p}\frac{1}{1-p^{-2}}=\zeta(2)\prod_{p|n} \biggl(1+\frac{1}{p}\biggr)=\frac{\pi^2}{6} \prod_{p|n} \biggl(1+\frac{1}{p}\biggr)\\
&= \frac{\pi^2}{6} \prod_{p|n:\,\, p\leq (\ln n)^{\alpha}} \biggl(1+\frac{1}{p}\biggr)\prod_{p|n:\,\, (\ln n)^{\alpha}< p\leq \ln n} \biggl(1+\frac{1}{p}\biggr)\prod_{p|n:\,\, p> \ln n} \biggl(1+\frac{1}{p}\biggr).
\end{align*} By Lemma \ref{L2} the last product does not exceed 5. It is well-known that (see, for example, \cite[Chapter 1]{Prachar})
\begin{equation}\label{B_2_ineq}
B_1 \ln x \leq \prod_{p\leq x} \biggl(1+\frac{1}{p}\biggr)\leq B_2\ln x,\qquad x\geq 2,
\end{equation}where $B_1 >0$ and $B_2>0$ are absolute constants. We obtain
\begin{align*}
&\prod_{p|n:\,\, (\ln n)^{\alpha}< p\leq \ln n} \biggl(1+\frac{1}{p}\biggr)\leq \prod_{(\ln n)^{\alpha}< p\leq \ln n} \biggl(1+\frac{1}{p}\biggr)\\
&\ \ \ \ \ =\prod_{p\leq \ln n} \biggl(1+\frac{1}{p}\biggr) \bigg/ \prod_{ p\leq (\ln n)^{\alpha}} \biggl(1+\frac{1}{p}\biggr)\leq \frac{B_2 \ln\ln n}{B_1\ln (\ln n)^{\alpha}}=\frac{B_2 \ln\ln n}{B_1 \alpha \ln\ln n}=\frac{B}{\alpha}.
\end{align*}Hence,
\[
\frac{n}{\varphi(n)}\leq \frac{5\pi^2B}{6\alpha}\prod_{p|n:\,\, p\leq (\ln n)^{\alpha}} \biggl(1+\frac{1}{p}\biggr)=
C(\alpha)\prod_{p|n:\,\, p\leq (\ln n)^{\alpha}} \biggl(1+\frac{1}{p}\biggr).
\]Lemma \ref{L3} is proved.

Let us continue the proof of Theorem \ref{T1}. We may assume that $M\geq \textup{exp} (2^{1/\alpha})$. We put
\[
y=(\ln M)^{\alpha}\quad\quad\text{and}\quad\quad S=\sum_{n=1}^{N}\left(\frac{a_n}{\varphi (a_n)}\right)^{s}.
\]We observe that $y\geq 2$. Let $1\leq n \leq N$. By Lemma \ref{L3}, we have
 \begin{align*}
 \frac{a_n}{\varphi (a_n)}&\leq C(\alpha) \prod_{p|a_n:\,\, p\leq (\ln a_n)^{\alpha}} \biggl(1+\frac{1}{p}\biggr)\\
 &\leq C(\alpha) \prod_{p|a_n:\,\, p\leq y} \biggl(1+\frac{1}{p}\biggr)= C(\alpha)\sum_{\substack{d|a_n:\\d\in \mathcal{M},\,P^{+}(d)\leq y}}\frac{1}{d}.
  \end{align*}Hence,
  \begin{align*}
  S&\leq (C(\alpha))^{s}\sum_{1\leq n \leq N}\biggl(\sum_{\substack{d_1|a_n:\\d_{1}\in \mathcal{M},\,P^{+}(d_{1})\leq y}}\frac{1}{d_1}\biggr)\cdots
  \biggl(\sum_{\substack{d_s|a_n:\\ d_{s}\in \mathcal{M},\,P^{+}(d_{s})\leq y}}\frac{1}{d_s}\biggr)\\
  &=(C(\alpha))^{s}\sum_{1\leq n \leq N}\sum_{\substack{d_1|a_n,\ldots,\,d_s|a_n\\
  d_1,\ldots,\,d_s \in \mathcal{M}\\
  P^{+}(d_1)\leq y,\ldots,\,P^{+}(d_s)\leq y}}\frac{1}{d_1\cdots  d_s}\\
  &\leq (C(\alpha))^{s}\sum_{\substack{1\leq d_1 \leq M\\
  P^{+}(d_1)\leq y\\
  d_1\in \mathcal{M}}}\dots \sum_{\substack{1\leq d_s \leq M\\
  P^{+}(d_s)\leq y\\
  d_s\in \mathcal{M}}} \frac{1}{d_1\cdots  d_s}\sum_{\substack{1\leq n \leq N:\\
  d_1|a_n,\ldots,\,d_s|a_n}} 1\\
  &=(C(\alpha))^{s}\sum_{\substack{1\leq d_1 \leq M\\
  P^{+}(d_1)\leq y\\
  d_1\in \mathcal{M}}}\dots \sum_{\substack{1\leq d_s \leq M\\
  P^{+}(d_s)\leq y\\
  d_s\in \mathcal{M}}} \frac{\omega([d_1,\ldots, d_s])}{d_1\cdots  d_s}.
  \end{align*} We observe that if $d$ and $d^{\prime}$ are positive integers and $d^{\prime}|d$, then
  \[
  \omega(d) \leq \omega (d^{\prime}).
  \]Let $d_1,\ldots, d_s$ be integers such that $1\leq d_i \leq M$, $P^{+}(d_i)\leq y$, $d_i\in \mathcal{M}$ for all $1\leq i \leq s$. Since
  \[
  P^{+}([d_1,\ldots, d_s])\,|\, [d_1,\ldots, d_s],
  \]we have
  \[
  \omega([d_1,\ldots, d_s]) \leq \omega\bigl(P^{+}([d_1,\ldots, d_s])\bigr).
  \]Since
  \[
  P^{+}([d_1,\ldots, d_s])=P^{+}(d_1\cdots d_s),
  \]we obtain
  \begin{align}
  S&\leq(C(\alpha))^{s}\sum_{\substack{1\leq d_1 \leq M\\
  P^{+}(d_1)\leq y\\
  d_1\in \mathcal{M}}}\dots \sum_{\substack{1\leq d_s \leq M\\
  P^{+}(d_s)\leq y\\
  d_s\in \mathcal{M}}} \frac{\omega(P^{+}(d_1\cdots d_s))}{d_1\cdots  d_s}\notag\\
  &\leq (C(\alpha))^{s}\sum_{\substack{d_{1}\in \mathcal{M}:\\P^{+}(d_1)\leq y}}\cdots \sum_{\substack{d_s\in \mathcal{M}:\\P^{+}(d_s)\leq y}}
  \frac{\omega(P^{+}(d_1\cdots d_s))}{d_1\cdots  d_s}= (C(\alpha))^{s} S^{\prime}.\label{T1:S_estimate}
  \end{align}

  It is easy to see that
  \begin{equation}\label{T1:S_prime}
  S^{\prime}=\omega(1)+\sum_{p\leq y} \omega(p) S_{p},
  \end{equation}where
  \[
  S_{p}=\sum_{\substack{d_{1},\ldots, d_{s}\in \mathcal{M}:\\
   P^{+}(d_{1})\leq p,\ldots, P^{+}(d_{s})\leq p,\\
   \text{and }\exists \tau:\ p|d_{\tau}}} \frac{1}{d_{1}\cdots d_{s}}.
  \] For a prime $p$ with $p\leq y$ and an integer $\tau$ with $1\leq \tau \leq s$, we put
  \[
  S_{p}(\tau)=\sum_{\substack{d_{1},\ldots, d_{s}\in \mathcal{M}:\\
   P^{+}(d_{1})\leq p,\ldots, P^{+}(d_{s})\leq p,\\
   \text{and } p|d_{\tau}}} \frac{1}{d_{1}\cdots d_{s}}.
  \]Applying \eqref{B_2_ineq}, we have (the product here is over prime $q$)
  \begin{align*}
  S_{p}(\tau)&\leq \frac{1}{p}\sum_{\substack{d_{1},\ldots, d_{s}\in \mathcal{M}:\\
   P^{+}(d_{1})\leq p,\ldots, P^{+}(d_{s})\leq p}} \frac{1}{d_{1}\cdots d_{s}}=
   \frac{1}{p} \Biggl(\sum_{\substack{d\in \mathcal{M}:\\
    P^{+}(d)\leq p}}\frac{1}{d}\Biggr)^{s}\\
    &=\frac{1}{p}\biggl(\prod_{q\leq p}\biggl(1+\frac{1}{q}\biggr)\biggr)^{s}\leq
    \frac{(B_{2}\ln p)^{s}}{p}.
  \end{align*}It is easy to see that
  \[
  S_{p}\leq \sum_{\tau=1}^{s}S_{p}(\tau).
  \]We obtain
  \[
  S_{p} \leq s \frac{(B_{2}\ln p)^{s}}{p}\leq \frac{(2B_{2}\ln p)^{s}}{p}.
  \] We may assume that $2B_{2}\geq 1$. Applying \eqref{T1:S_prime} and taking into account that $\omega(1)=N$ and $y=(\ln M)^{\alpha}$, we obtain
  \[
  S^{\prime} \leq (2B_{2})^{s}\biggl(N+\sum_{p\leq (\ln M)^{\alpha}}\frac{\omega(p)(\ln p)^{s}}{p}\biggr).
  \]By \eqref{T1:S_estimate}, we have
  \begin{gather*}
  S\leq \bigl(2B_{2}C(\alpha)\bigr)^{s}\biggl(N+\sum_{p\leq (\ln M)^{\alpha}}\frac{\omega(p)(\ln p)^{s}}{p}\biggr)\\
  =\bigl( \widetilde{C}(\alpha)\bigr)^{s} \biggl(N+\sum_{p\leq (\ln M)^{\alpha}}\frac{\omega(p)(\ln p)^{s}}{p}\biggr),
  \end{gather*}where $\widetilde{C}(\alpha)>0$ is a constant, depending only on $\alpha$. Theorem \ref{T1} is proved.\\[-5pt]

  \textsc{Proof of Theorem \ref{T2}.} We will choose the constant $M_0> 0$, depending on the sequence $\alpha(M)$, later; it will be large enough. For now, let $M_0$ satisfy the following conditions: $M_0 \geq 100$, and
  \[
  \alpha(M)\leq \frac{1}{2},\quad\quad  (\ln M)^{\alpha(M)}\leq \frac{\ln M}{4}
  \]for any $M\geq M_0$. Let $M\geq M_0$. We put
  \[
  y= (\ln M)^{\alpha(M)},\quad\quad z=\frac{\ln M}{2}.
  \]Then $2\leq y \leq z/2.$ We define
  \[
  A=\{1\leq n \leq M: \text{$p|n$ for any $p\in (y,z]$, and $p\nmid n$ for any $p\leq y$}\}.
  \]We put
  \[
  Q=\prod_{y<p\leq z} p.
  \] We have
   \[
   \ln Q=\sum_{y<p\leq z} \ln p= \theta(z) - \theta(y),\ \text{where }\theta(x)=\sum_{p\leq x}\ln p.
   \] Since (see, for example, \cite[Chapter 3]{Prachar})
   \begin{equation}\label{teta_equiv_x}
   \lim_{x\to +\infty}\frac{\theta(x)}{x}=1,
   \end{equation}we see that there is an absolute constant $c_1>0$ such that $\theta (x) \geq x/2$ for all $x\geq c_1$. We may assume that $M_0 > \textup{exp} (2 c_1)$; therefore $z=(\ln M)/2\geq c_1$, and, hence, $\theta (z)\geq z/2 = (\ln M)/4$. It follows from \eqref{teta_equiv_x} that $ \theta(x)\leq b x$ for all $x\geq 2$, where $b>0$ is an absolute constant. Since $y\geq 2$, we obtain
   \[
   \theta (y) \leq b y = b (\ln M)^{\alpha (M)}\leq b (\ln M)^{1/2}.
   \]Hence,
   \[
   \ln Q\geq \frac{\ln M}{4} - b (\ln M)^{1/2}\geq 100,
   \]if $M_0$ is chosen large enough. In particular, we see that $\Omega = \{p:  y< p \leq z\}\neq \emptyset$ and $Q\geq 100$.
   It follows from \eqref{teta_equiv_x} that there is an absolute constant $c_2>0$ such that $\theta (x) \leq (3/2)\,x$ for all $x\geq c_2$. We may assume that $M_0 > \textup{exp}(2c_2)$, and, hence, $z=(\ln M)/2 \geq c_2$.
   Therefore
   \[
   \ln Q\leq \theta (z)\leq \frac{3}{2}\,z=\frac{3}{4}\,\ln M=\ln (M^{3/4}).
   \]We obtain
   \[
   Q\leq M^{3/4}<M.
   \]We see that $Q\in A$, and hence $A\neq\emptyset$.

   Thus, $A\subset\{1,\ldots, M\}$, $A\neq \emptyset$, and
   \[
   \#\{n\in A: n\equiv 0\ (\textup{mod } p)\}=0
   \]for any prime $p\leq y=(\ln M)^{\alpha(M)}$. If $n\in A$, then $n>1$, since $Q|n$ and $Q\geq100$.

    Let $n\in A$. Then
    \[
    \frac{n}{\varphi (n)}=\prod_{p|n}\biggl(1-\frac{1}{p}\biggr)^{-1}\geq\prod_{y<p\leq z}\biggl(1-\frac{1}{p}\biggr)^{-1}=\prod_{p\leq z}\biggl(1-\frac{1}{p}\biggr)^{-1}\bigg/ \prod_{p\leq y}\biggl(1-\frac{1}{p}\biggr)^{-1}.
    \]Since
    \[
    D_1 \ln x\leq \prod_{p \leq x} \biggl(1-\frac{1}{p}\biggr)^{-1}\leq D_2 \ln x,\qquad x\geq 2,
    \]where $D_1 >0$ and $D_2 >0$ are absolute constants, and $2 \leq y \leq z/2< z$, we have
    \begin{align*}
    &\prod_{p\leq y}\biggl(1-\frac{1}{p}\biggr)^{-1} \leq D_2 \ln y= D_2\alpha(M)\ln\ln M,\\
    &\prod_{p\leq z}\biggl(1-\frac{1}{p}\biggr)^{-1}\geq D_1 \ln z= D_1 \ln \biggl(\frac{1}{2}\ln M\biggr)\geq \frac{D_1}{2}\ln\ln M.
    \end{align*}Hence,
    \[
    \frac{n}{\varphi(n)}\geq \frac{(D_1 / 2)\ln\ln M}{D_2 \alpha(M)\ln\ln M}=\frac{c}{\alpha(M)}.
    \]We obtain
    \[
    \sum_{n\in A}\frac{n}{\varphi (n)}\geq \frac{c}{\alpha (M)} \#A,
    \]where $c>0$ is an absolute constant. Theorem \ref{T2} is proved.\\[-5pt]

    \textsc{Proof of Theorem \ref{T3}.} We have
\begin{align*}
R(n)=a_k n^k+\ldots+ a_0,\qquad a_k\neq 0,\quad \delta=(a_0,\ldots,a_k),\\
 |a_i|\leq x,\ i=0,\ldots, k,\qquad (\ln x)^{\varepsilon} \leq z \leq x,\quad x\geq2.
\end{align*} We put
\[
\Omega=\{-z\leq n \leq z: R(n)\neq 0\}.
\] Suppose that $\Omega \neq \emptyset$. We are going to obtain an upper bound for
\[
S =\sum_{n\in \Omega} \biggl(\frac{|R(n)|}{\varphi (|R(n)|)}\biggr)^{s}.
\]We put
\[
\widetilde{R}(n)=\frac{1}{\delta} R(n).
\]Then
\begin{align*}
\widetilde{R}(n)=\widetilde{a}_{k}n^k&+\ldots+\widetilde{a}_0,\qquad (\widetilde{a}_{0},\ldots, \widetilde{a}_k) = 1,\\
&|\widetilde{a}_{i}|\leq x,\qquad i=0,\ldots, k.
\end{align*}Since $\varphi (mn)\geq \varphi(m)\varphi(n)$
 for all positive integers $m$ and $n$, we obtain
\begin{equation}\label{T3:S.S.wave.DEPENDENCE}
S= \sum_{n\in \Omega} \biggl(\frac{\delta|\widetilde{R}(n)|}{\varphi (\delta|\widetilde{R}(n)|)}\biggr)^{s}\leq \biggl(\frac{\delta}{\varphi(\delta)}\biggr)^{s}\sum_{n\in \Omega} \biggl(\frac{|\widetilde{R}(n)|}{\varphi (|\widetilde{R}(n)|)}\biggr)^{s}= \biggl(\frac{\delta}{\varphi(\delta)}\biggr)^{s} \widetilde{S}.
\end{equation} Let $n\in \Omega$. Since $|\widetilde{a}_i|\leq x$, $|n|\leq z\leq x$, and $x\geq 2$, we have
\begin{align*}
|\widetilde{R}(n)|&=|\widetilde{a}_k n^k+\ldots+\widetilde{a}_0|\leq x(1+x+\ldots+ x^k)\\
&=x\,\frac{x^{k+1}-1}{x-1}\leq x\,\frac{x^{k+1}}{x/2}= 2 x^{k+1}\leq x^{k+2}\leq x^{3k}.
\end{align*} We put $M= x^{3k}$.
 We have proved that if $n\in\Omega$, then $|\widetilde{R}(n)| \leq M$.

  Let $c(\varepsilon)>0$ be such that the following holds:

\textup{i}) $c(\varepsilon)\geq 30$;

\textup{ii}) $\ln x \geq 2^{4/\varepsilon}$ for $x\geq c(\varepsilon)$;

\textup{iii}) $\bigl(3 (\ln x)^2\bigr)^{\varepsilon/4} \leq (\ln x)^{\varepsilon}$ for $x\geq c(\varepsilon)$.

Let $x \geq c(\varepsilon)$. Let us consider two cases.

1) Let $k\geq \ln x$. Let $n\in \Omega$. Then
 \begin{align*}
 1&\leq \frac{|\widetilde{R}(n)|}{\varphi(|\widetilde{R}(n)|)} \leq c\ln\ln(|\widetilde{R}(n)|+2)\leq
 c\ln\ln(x^{3k}+2)\leq c\ln\ln(x^{4k})\\
 &=c(\ln k+ \ln\ln x+ \ln 4)\leq c(\ln k+ 3\ln\ln x)\leq c(4 \ln k)=c_1 \ln k.
 \end{align*}Hence,
 \[
 \widetilde{S}=\sum_{n\in \Omega} \biggl(\frac{|\widetilde{R}(n)|}{\varphi(|\widetilde{R}(n)|)}\biggr)^s \leq (c_1 \ln k)^{s} \# \Omega.
 \] We recall that  $(\ln x)^{\varepsilon} \leq z \leq x$, and, hence, $z\geq 1$. Since
   \begin{align*}
 \# \Omega = \#\{ -z\leq n \leq z: \widetilde{R}(n)\neq 0\}\leq 2z+1 \leq 3z,
 \end{align*}we have $\widetilde{S}\leq (c_2 \ln k)^s z$. Since $S\leq (\delta/\varphi(\delta))^{s}\widetilde{S}$, we obtain
   \[
  S=\sum_{n\in \Omega} \biggl(\frac{|R(n)|}{\varphi(|R(n)|)}\biggr)^s \leq \Bigl(c_2\,\frac{\delta}{\varphi(\delta)}\ln k\Bigr)^{s} z.
  \]

  2) Let $k< \ln x$. We have
    \[
  \ln M = 3 k\ln x\leq 3(\ln x)^{2}.
  \]We take $\alpha = \varepsilon/4$.
  Using conditions on $c(\varepsilon)$, we obtain
  \begin{gather*}
  \ln M =3 k\ln x\geq 3 \ln x \geq 2^{1/\alpha},\\
  2\leq (\ln M)^{\alpha}\leq \bigl(3(\ln x)^2\bigr)^{\alpha}= \bigl(3(\ln x)^2\bigr)^{\varepsilon / 4}\leq(\ln x)^{\varepsilon}\leq z.
  \end{gather*}We define
  \[
  \omega(d)=\#\{n\in \Omega: \widetilde{R}(n)\equiv 0\ (\textup{mod }d)\}
  \] for any positive integer $d$. Let $p$ be a prime number. We have
  \[
  \omega(p)\leq \#\{-z\leq n\leq z: \widetilde{R}(n)\equiv 0\ (\textup{mod }p)\}.
  \]Since $(\widetilde{a}_0,\ldots, \widetilde{a}_k)=1$, there is a number $\widetilde{a}_i$ such that $\widetilde{a}_i\not \equiv 0$ (mod $p$). Therefore the number of solutions of the congruence
  \begin{equation}\label{R=0}
  \widetilde{R}(n)\equiv 0\quad (\textup{mod } p)
  \end{equation}does not exceed $k$. It is clear that the number of solutions does not exceed $p$. Thus, the number of solutions of the congruence \eqref{R=0} does not exceed $\min (p,k)$. Let $m_1<\ldots< m_t$ be all numbers from $\{1,\ldots, p\}$ satisfying the congruence \eqref{R=0} (hence, $t\leq \min (p,k)$). Let $1 \leq j \leq t$. We have
   \[
  \#\{-z\leq n\leq z: n\equiv m_j\ (\textup{mod\ }p)\}\leq \frac{2z}{p}+1.
  \]We obtain
  \[
  \omega(p)\leq \min(p,k)\biggl(\frac{2z}{p}+1\biggr).
  \]If $p\leq (\ln M)^{\alpha}$, then $p\leq z$. Therefore
  \begin{equation}\label{w(p)}
  \omega(p)\leq \min (p, k)\,\frac{3z}{p}
  \end{equation} for any $p\leq(\ln M)^{\alpha}$. Applying Theorem \ref{T1}, we obtain
  \[
  \widetilde{S}=\sum_{n\in \Omega} \biggl(\frac{|\widetilde{R}(n)|}{\varphi(|\widetilde{R}(n)|)}\biggr)^s\leq (C(\alpha))^{s}
  \biggl(\#\Omega +\sum_{p\leq (\ln M)^{\alpha}}\frac{\omega(p) (\ln p)^{s}}{p}\biggr).
  \]Since $\alpha=\varepsilon /4$, we see that $C(\alpha)=C(\varepsilon)>0$ is a constant depending only on $\varepsilon$. Applying \eqref{w(p)} and taking into account that $\#\Omega \leq 2z+1 \leq 3z$, we obtain
  \begin{align}
   \widetilde{S}&=\sum_{n\in \Omega} \biggl(\frac{|\widetilde{R}(n)|}{\varphi(|\widetilde{R}(n)|)}\biggr)^s \leq (C(\varepsilon))^{s}
  \biggl(3z +\sum_{p\leq (\ln M)^{\alpha}}\frac{\min(p,k) 3z (\ln p)^{s}}{p^2}\biggr)\notag\\
  &\leq (C(\varepsilon))^{s}
  \biggl(3z +\sum_{p}\frac{\min(p,k) 3z (\ln p)^{s}}{p^2}\biggr)= (C(\varepsilon))^{s} 3z
  \biggl(1 +\sum_{p}\frac{\min(p,k) (\ln p)^{s}}{p^2}\biggr).\label{T3.Est.S.Wave}
  \end{align}

  \begin{lemma}\label{L:Gamma}
  Let $s$ be a positive integer, $x$ be a real number with $x\geq 1$. Then
  \begin{equation}\label{GAMMA.INEQ}
  \int_{x}^{+\infty}t^{s-1}e^{-t}\,dt\leq s!\,x^{s-1}e^{-x}.
  \end{equation}
  \end{lemma}

  \textsc{Proof of Lemma \ref{L:Gamma}.} We put
  \[
  \Gamma(s,x)=\int_{x}^{+\infty}t^{s-1}e^{-t}\,dt.
  \]We have
  \begin{align}
  \Gamma(s,x)&=\int_{x}^{+\infty}t^{s-1}e^{-t}\,dt=-\int_{x}^{+\infty}t^{s-1}\,d e^{-t}\notag\\
  &=-t^{s-1}e^{-t}\Bigr|_{x}^{+\infty}+(s-1)\int_{x}^{+\infty}e^{-t}t^{s-2}\,dt=x^{s-1}e^{-x}+(s-1)\Gamma(s-1,x).\label{L:GAMMA_INDUCTION}
  \end{align}

  Now we induct on $s$. If $s=1$, we have $\Gamma(1,x)=e^{-x}$ and \eqref{GAMMA.INEQ} holds. Suppose that $s\geq 2$ and the claim is true for $s-1$. Then \eqref{L:GAMMA_INDUCTION} and the induction hypothesis give us
  \begin{align*}
  \Gamma(s,x)&=x^{s-1}e^{-x}+(s-1)\Gamma(s-1,x)\leq x^{s-1}e^{-x} + (s-1)(s-1)!\,x^{s-2}e^{-x}\\
  &=s!\,x^{s-1}e^{-x}\left(\frac{1}{s!}+ \frac{1-1/s}{x}\right)\leq s!\,x^{s-1}e^{-x} \left(1 -\frac{1}{s}+\frac{1}{s!} \right) \leq s!\,x^{s-1}e^{-x},
  \end{align*}since $x\geq 1$. The claim follows. Lemma \ref{L:Gamma} is proved.

  \begin{lemma}\label{L:primes}
  Let $k$ and $s$ be positive integers. Then
  \begin{align}
  &\sum_{p\leq k}\frac{(\ln p)^{s}}{p} \leq c (\ln k)^s,\label{BASIC.PRIMES.1}\\
  &\sum_{p>k}\frac{(\ln p)^{s}}{p^2} \leq c s!\,\frac{(\ln (k+2))^{s-1}}{k},\label{BASIC.PRIMES.2}
  \end{align}where $c$ is a positive absolute constant.
  \end{lemma}

  \textsc{Proof of Lemma \ref{L:primes}.} For any positive integer $n$, we put
  \[
  a_{n}= \begin{cases}
  (\ln n)^{s}    &\text{if $n\in \mathbb{P}$;}\\
   0             &\text{otherwise}.
          \end{cases}
  \]We set
  \[
  A(x)=\sum_{n\leq x} a_{n}.
  \]For any real number $x\geq 2$, we have
  \[
  A(x)=\sum_{p\leq x} (\ln p)^{s}\leq (\ln x)^{s}\pi(x)\leq (\ln x)^{s} c\,\frac{x}{\ln x}=
  c x (\ln x)^{s-1},
  \]where $c$ is a positive absolute constant. Since $A(x)=0$ for $1\leq x < 2$, we obtain $A(x) \leq c x (\ln x)^{s-1}$ for $x\geq 1$.

  I) We first prove \eqref{BASIC.PRIMES.1}. We denote the sum in \eqref{BASIC.PRIMES.1} by $S_1$. Suppose that $k \geq 2$. We put $f(x)=1/x$. Applying partial summation (see, for example, \cite[Theorem 2.1.1]{Murty}), we have
  \[
  S_{1}=\sum_{n\leq k} a_{n}f(n)= A(k)f(k) - \int_{1}^{k} A(x)f'(x)\, dx.
  \] We have
  \[
  A(k)f(k)\leq c (\ln k)^{s-1}\leq c\, \frac{\ln k}{\ln 2}(\ln k)^{s-1} \leq 2c (\ln k)^{s}
  \] and
  \begin{align*}
  -\int_{1}^{k}A(x)f'(x)\, dx&= \int_{1}^{k}\frac{A(x)}{x^{2}}\,dx\leq
  c\int_{1}^{k} \frac{(\ln x)^{s-1}}{x}\, dx\\
  &=c\int_{0}^{\ln k} t^{s-1}\, dt= \frac{c}{s} (\ln k)^{s}\leq c (\ln k)^{s}.
  \end{align*} We obtain $S_{1} \leq (3c) (\ln k)^{s}$.

  Finally, if $k=1$, then $S_{1}= 0$, and the previous inequality also holds. The inequality \eqref{BASIC.PRIMES.1} is proved.

  II) Now we prove \eqref{BASIC.PRIMES.2}. We denote the sum in \eqref{BASIC.PRIMES.2} by $S_{2}$. We put $f(x)=1/x^{2}$. Applying partial summation, we have
  \[
  \sum_{n\leq u} a_{n}f(n) = A(u)f(u) - \int_{1}^{u}A(x)f'(x)\, dx
  \]for any real number $u\geq 1$. Since $A(u)f(u)\to 0$ as $u\to +\infty$, we obtain
  \[
  \sum_{n=1}^{+\infty}a_{n}f(n)= - \int_{1}^{+\infty} A(x)f'(x)\, dx.
  \]Also we have
  \[
  \sum_{n\leq k}a_{n}f(n)= A(k)f(k) - \int_{1}^{k}A(x)f'(x)\, dx.
  \]Therefore
  \begin{align*}
  S_{2}&=\sum_{n\geq k+1} a_{n}f(n)= - \int_{k}^{+\infty} A(x)f'(x)\, dx - A(k)f(k)\leq- \int_{k}^{+\infty} A(x)f'(x)\, dx\\
  &= 2\int_{k}^{+\infty}\frac{A(x)}{x^{3}}\, dx\leq
  2c\int_{k}^{+\infty}\frac{(\ln x)^{s-1}}{x^{2}}\,dx= 2 c \int_{\ln k}^{+\infty} t^{s-1}e^{-t}\, dt= 2c I_{k}.
  \end{align*} If $k\geq 3$, then by Lemma \ref{L:Gamma} we have
  \[
  I_{k}\leq s!\,(\ln k)^{s-1} e^{-\ln k}=s!\,\frac{(\ln k)^{s-1}}{k}\leq s!\,\frac{(\ln (k+2))^{s-1}}{k},
  \] and
  \[
  I_{k} \leq \int_{0}^{+\infty} t^{s-1}e^{-t}\, dt= \Gamma(s)= (s-1)!\,\leq s!\,2 \frac{(\ln (k+2))^{s-1}}{k},
  \]if $k\in \{1, 2\}$. We see that
  \[
  I_{k} \leq s!\,2 \frac{(\ln (k+2))^{s-1}}{k}
  \]for any positive integer $k$. We obtain
  \[
  S_{2} \leq 4c s!\,\frac{(\ln (k+2))^{s-1}}{k},
  \]and the inequality \eqref{BASIC.PRIMES.2} is proved. Lemma \ref{L:primes} is proved.

  We may assume that $c\geq 1$, where $c$ is the constant in Lemma \ref{L:primes}. Applying Lemma \ref{L:primes} and taking into account that $\ln (k+2) \leq 2 \ln (k+1)$, we obtain
  \begin{align*}
  \sum_{p\leq k}\frac{(\ln p)^{s}}{p} &\leq c \bigl(\ln k\bigr)^s \leq c \bigl(\ln (k+1)\bigr)^s \leq
  c^{s} \bigl(\ln (k+1)\bigr)^s \leq c^{s} s!\,\bigl(\ln (k+1)\bigr)^s,\\
  k\sum_{p>k}\frac{(\ln p)^{s}}{p^2} &\leq c s!\,\bigl(\ln (k+2)\bigr)^{s-1}\leq c s!\,\bigl(\ln (k+2)\bigr)^{s}\leq
  c s!\,2^{s} \bigl(\ln (k+1)\bigr)^{s}\\
  &\leq (2c)^{s} s!\,\bigl(\ln (k+1)\bigr)^{s}.
  \end{align*}Since $c^{s}+(2c)^{s}\leq 2 (2c)^{s}\leq (4c)^{s}$, we have
  \[
  \sum_{p}\frac{\min(p,k) (\ln p)^{s}}{p^2}= \sum_{p\leq k}\frac{(\ln p)^{s}}{p}+k\sum_{p>k}\frac{(\ln p)^{s}}{p^2}
  \leq (4c)^{s} s!\,\bigl(\ln (k+1)\bigr)^{s}.
  \] Since
  \[
  1\leq \bigl(\ln (k+2)\bigr)^{s}\leq 2^{s}\bigl(\ln (k+1)\bigr)^{s}\leq (2c)^{s}s!\,\bigl(\ln (k+1)\bigr)^{s},
  \]we obtain
  \begin{equation}\label{T3:BASIC.SUM}
  1+ \sum_{p}\frac{\min(p,k) (\ln p)^{s}}{p^2}\leq (8c)^{s}s!\,\bigl(\ln (k+1)\bigr)^{s}.
  \end{equation}By \eqref{T3.Est.S.Wave}, we have
  \begin{align*}
  \widetilde{S}&\leq (C(\varepsilon))^{s} 3z (8c)^{s}s!\,\bigl(\ln (k+1)\bigr)^{s}\leq
  (24c C(\varepsilon))^{s} z s!\,\bigl(\ln (k+1)\bigr)^{s}\\
  &=(C_{1}(\varepsilon))^{s} z s!\,\bigl(\ln (k+1)\bigr)^{s},
  \end{align*}where $C_{1}(\varepsilon)>0$ is a constant, depending only on $\varepsilon$. By \eqref{T3:S.S.wave.DEPENDENCE}, we obtain
  \[
  S=\sum_{n\in \Omega} \biggl(\frac{|R(n)|}{\varphi(|R(n)|)}\biggr)^s\leq \biggl( C_1(\varepsilon)\frac{\delta}{\varphi(\delta)} \ln (k+1)\biggr)^{s}s!\,z.
  \]Thus, Theorem \ref{T3} in the case $x\geq c(\varepsilon)$ is proved. Since $n/\varphi(n) \leq c \ln\ln (n+2)$, the claim in the case $3 \leq x< c(\varepsilon)$ is trivial. Theorem \ref{T3} is proved.\\[-5pt]

  \textsc{Proof of Corollary \ref{Cor1}.} We put
  \[
  x_{0}(R)=\max (|a_{0}|,\ldots, |a_{k}|)+10.
  \]We see that $x_{0}(R)$ is a positive constant depending only on $R$. Suppose that $x\geq x_{0}(R)$. Applying Theorem \ref{T3} with $\varepsilon=1/2$ and $z=x$, we obtain
  \[
  \sum_{\substack{-x \leq n \leq x\\ R(n)\neq 0}} \left(\frac{|R(n)|}{\varphi(|R(n)|)}\right)^{s}\leq
  \left(C(1/2)\frac{\delta}{\varphi(\delta)}\ln (k+1)\right)^{s}s!\,x=
  (c_{1}(R))^{s}s!\,x,
  \]where $\delta=(a_{0},\ldots, a_k)$ and $c_{1}(R)= C(1/2)(\delta/\varphi(\delta))\ln (k+1)$ is a positive constant depending only on $R$.

  Suppose that $1\leq x < x_{0}(R)$. Let $\Omega=\{n: -x \leq n \leq x\text{ and $R(n)\neq0$}\}\neq\emptyset$. We put
  \[
  m(R)=\max_{-x_{0}(R)\leq n \leq x_{0}(R)} |R(n)|+10.
  \]We see that $m(R)$ is a positive constant depending only on $R$. For any integer $n$ such that $-x_{0}(R)\leq n \leq x_{0}(R)$ and $R(n)\neq 0$, we have
  \[
  \frac{|R(n)|}{\varphi(|R(n)|)} \leq |R(n)|\leq m(R).
  \]We obtain
  \begin{align*}
  S&=\sum_{\substack{-x \leq n \leq x\\ R(n)\neq 0}} \left(\frac{|R(n)|}{\varphi(|R(n)|)}\right)^{s}\leq
  \sum_{\substack{-x_{0}(R) \leq n \leq x_{0}(R)\\ R(n)\neq 0}} \left(\frac{|R(n)|}{\varphi(|R(n)|)}\right)^{s}\\
  &\leq  (m(R))^{s} (2 x_{0}(R)+1)\leq
  \bigl(3x_{0}(R)m(R)\bigr)^{s} s!\,x= (c_{2}(R))^{s} s!\,x,
  \end{align*}where $c_{2}(R)= 3x_{0}(R)m(R)$ is a positive constant depending only on $R$. If $\Omega = \emptyset$, then $S=0$. Thus, the claim follows with $C(R)=\max (c_{1}(R), c_{2}(R))$. Corollary \ref{Cor1} is proved.\\[-5pt]

  \textsc{Proof of Theorem \ref{T4}.} We assume that $x\geq c(\varepsilon)$, where $c(\varepsilon)$ is a positive constant depending only on $\varepsilon$; we will choose the constant $c(\varepsilon)$ later, and it will be large enough.

   Suppose that $\Omega:=\{-z \leq b\leq z: \Delta_{L}\neq 0\}\neq \emptyset$. We put $\Delta (b):= \prod_{i=1}^{k} |b_i - b|$. We have
   \[
   \frac{\Delta_{L}}{\varphi (\Delta_{L})}\leq \frac{a^{k+1} \Delta (b)}{\varphi(a^{k+1})\varphi (\Delta(b))}= \frac{a}{\varphi(a)} \frac{\Delta(b)}{\varphi(\Delta(b))}.
   \]Therefore
   \begin{equation}\label{T4.Sum.DELTA.L}
   \sum_{b\in \Omega} \biggl(\frac{\Delta_{L}}{\varphi (\Delta_{L})}\biggr)^{s} \leq \biggl(\frac{a}{\varphi(a)}\biggr)^{s}  \sum_{b\in \Omega} \biggl(\frac{\Delta (b)}{\varphi (\Delta (b))}\biggr)^{s}.
   \end{equation}To obtain an upper bound for the last sum in \eqref{T4.Sum.DELTA.L} we are going to apply Theorem \ref{T1} with $\alpha=\varepsilon/4$.

    We put
    \[
    R(n)=(n-b_1)\cdots (n-b_k)=n^{k}+a_{k-1}n^{k-1}+\ldots+ a_0.
    \] Given a prime $p$, we have
    \[
    \omega(p)=\#\{b\in \Omega: \Delta(b)\equiv 0\text{ (mod $p$)}\}\leq \#\{-z \leq b \leq z: R(b)\equiv 0\text{ (mod $p$)}\}.
    \]The number of solutions of the congruence
    \begin{equation}\label{T4.R.CONGRUENCE}
    R(n)\equiv 0\ \text{ (mod $p$)}
    \end{equation}does not exceed $k$, and it trivially does not exceed $p$. We obtain the number of solutions of the congruence \eqref{T4.R.CONGRUENCE} does not exceed $\min (p,k)$. Let $m_1<\ldots < m_t$ be all numbers from $\{1,\ldots, p\}$ satisfying the congruence \eqref{T4.R.CONGRUENCE} (hence, $t\leq \min (p,k)$). Let $1 \leq j \leq t$. We have
    \[
    \#\{-z\leq b \leq z: b\equiv m_j\ \text{(mod $p$)}\}\leq \frac{2z}{p}+1.
    \]Hence,
    \begin{equation}\label{T4.omega.EST}
    \omega(p) \leq t \biggl(\frac{2z}{p}+1\biggr)\leq \min (p,k) \biggl(\frac{2z}{p}+1\biggr).
    \end{equation}

    Let $b\in \Omega$. Since
    \[
    |b-b_i|\leq |b|+|b_i|\leq z+x\leq 2x
    \] for all $1 \leq i \leq k$, we obtain $|R(b)| \leq (2x)^{k}=:M$.

    Let us consider two cases.

    1) Let $k\geq \ln x$. We may assume that $c(\varepsilon) \geq 30$, and therefore $x \geq 30$. We see that $(2x)^{l}+2\leq (3x)^{l}$ for any positive integer $l$, and $\ln\ln (3t)\leq c_0 \ln\ln t$ for any $t\geq 30$, where $c_0$ is a positive absolute constant. Let $b\in \Omega$. We have
    \begin{align*}
    \frac{\Delta(b)}{\varphi(\Delta(b))}&=\frac{|R(b)|}{\varphi(|R(b)|)}\leq c_{1}\ln\ln(|R(b)|+2)\leq
    c_{1}\ln\ln\bigl((2x)^{k}+2\bigr)\leq c_{1}\ln\ln\bigl((3x)^{k}\bigr)\\
    &=c_1\bigl(\ln k+ \ln\ln (3x)\bigr)\leq
    c_1(\ln k+ c_0 \ln\ln x)\leq c_1(\ln k+ c_0 \ln k)=c_2 \ln k,
    \end{align*} where $c_1$, $c_2= c_1 (1+c_0)$ are positive absolute constants. We observe that $z \geq (\ln x)^{\varepsilon}\geq 1$, and hence
    \[
    \#\Omega \leq 2z+1\leq 3z.
    \]We obtain
    \[
    \sum_{b\in \Omega} \biggl(\frac{\Delta(b)}{\varphi(\Delta (b))}\biggr)^{s} \leq (c_2 \ln k)^{s}\#\Omega \leq
    (c_2 \ln k)^{s} 3 z.
    \] From \eqref{T4.Sum.DELTA.L} we obtain
    \[
    \sum_{b\in \Omega} \biggl(\frac{\Delta_{L}}{\varphi(\Delta_{L})}\biggr)^{s}\leq \biggl(c  \frac{a}{\varphi (a)} \ln k \biggr)^{s} z,
    \]where $c$ is a positive absolute constant.

    2) Let $k < \ln x$. Let $c(\varepsilon)$ be such that $\ln (2t)\geq 2^{4/\varepsilon}$, $(\ln (2t)\ln t)^{\varepsilon/4}\leq (\ln t)^{\varepsilon}$ for any $t\geq c(\varepsilon)$. We have
    \begin{gather*}
    \ln M = k \ln (2x) \geq \ln (2x)\geq 2^{4/\varepsilon},\\
    (\ln M)^{\varepsilon/4}\leq \bigl(\ln (2x)\ln x\bigr)^{\varepsilon/4}\leq (\ln x)^{\varepsilon}\leq z.
    \end{gather*}We see that $2\leq (\ln M)^{\varepsilon/4}\leq z$, and hence $z/p \geq 1$ for any prime $p \leq (\ln M)^{\varepsilon/4}$. From \eqref{T4.omega.EST} we obtain
    \[
    \omega(p) \leq  \min (p,k) \frac{3z}{p}
    \]for any prime $p \leq (\ln M)^{\varepsilon/4}$. Applying Theorem \ref{T1}, we have
    \begin{align*}
    \sum_{b\in \Omega} \biggl(\frac{\Delta(b)}{\varphi(\Delta (b))}\biggr)^{s}&\leq (C(\varepsilon))^{s}
    \biggl(\#\Omega +\sum_{p\leq (\ln M)^{\varepsilon/4}} \frac{\omega(p)(\ln p)^{s}}{p}\biggr)\\
    &\leq (C(\varepsilon))^{s}
    \biggl(3z +\sum_{p\leq (\ln M)^{\varepsilon/4}} \frac{3z \min (p,k) (\ln p)^{s}}{p^2}\biggr)\\
    &\leq (C(\varepsilon))^{s} 3z
    \biggl(1 +\sum_{p} \frac{\min (p,k) (\ln p)^{s}}{p^2}\biggr).
    \end{align*}By \eqref{T3:BASIC.SUM}, we obtain
    \[
    \sum_{b\in \Omega} \biggl(\frac{\Delta(b)}{\varphi(\Delta (b))}\biggr)^{s}\leq \biggl(C_{1}(\varepsilon) \ln (k+1)\biggr)^{s} s!\,z,
    \]where $C_{1}(\varepsilon) >0$ is a constant, depending only on $\varepsilon$. From \eqref{T4.Sum.DELTA.L} we obtain
    \[
    \sum_{b\in \Omega} \biggl(\frac{\Delta_{L}}{\varphi(\Delta_{L})}\biggr)^{s}\leq
    \biggl(C_{1}(\varepsilon) \frac{a}{\varphi (a)} \ln (k+1)\biggr)^{s} s!\,z.
    \]Thus, Theorem \ref{T4} in the case $x\geq c(\varepsilon)$ is proved.

    Let $3 \leq x < c(\varepsilon)$. Given $b\in \Omega$, we have
    \[
    |b - b_i| \leq |b|+|b_i| \leq z + x\leq 2x \leq 2 c(\varepsilon)
    \] for all $1 \leq i \leq k$, and hence $\Delta (b)\leq \bigl(2 c(\varepsilon)\bigr)^{k}$. We obtain
    \begin{align*}
    \frac{\Delta (b)}{\varphi (\Delta (b))} &\leq c \ln\ln \bigl(\Delta (b) + 2\bigr)\leq
    c \ln\ln \bigl((2 c(\varepsilon))^{k} + 2\bigr)\leq c \ln\ln \bigl((3 c(\varepsilon))^{k}\bigr)\\
    &= c \bigl(\ln k+ \ln\ln (3 c(\varepsilon))\bigr)\leq c\bigl(\ln (k+1) + 2\ln\ln (3 c(\varepsilon)) \ln (k+1)\bigr)\\
    &= c_{1}(\varepsilon) \ln (k+1).
    \end{align*}  We obtain
    \[
    \sum_{b\in \Omega} \biggl(\frac{\Delta (b)}{\varphi (\Delta (b))}\biggr)^{s}\leq \bigl(c_{1}(\varepsilon) \ln (k+1)\bigr)^{s}\#\Omega .
    \] We have
    \[
    \# \Omega \leq 2z+1\leq 2x+1\leq 2 c(\varepsilon)+1\leq 3 c(\varepsilon).
    \]We observe that $z\geq 1$, since $z \geq (\ln x)^{\varepsilon}$ and $x \geq 3$. We obtain
    \begin{align*}
    \sum_{b\in \Omega} \biggl(\frac{\Delta (b)}{\varphi (\Delta (b))}\biggr)^{s}&\leq \bigl(c_{1}(\varepsilon) \ln (k+1)\bigr)^{s} 3 c(\varepsilon) z\leq
    \bigl(c_{1}(\varepsilon) 3 c(\varepsilon)  \ln (k+1)\bigr)^{s} z\\
    &= \bigl(c_{2}(\varepsilon) \ln (k+1)\bigr)^{s} z\leq \bigl(c_{2}(\varepsilon) \ln (k+1)\bigr)^{s} s!\,z,
    \end{align*}where $c_{2}(\varepsilon) >0$ is a constant depending only on $\varepsilon$. From \eqref{T4.Sum.DELTA.L} we obtain
    \[
    \sum_{b\in \Omega} \biggl(\frac{\Delta_{L}}{\varphi(\Delta_{L})}\biggr)^{s}\leq
    \biggl(c_{2}(\varepsilon) \frac{a}{\varphi (a)} \ln (k+1)\biggr)^{s} s!\,z.
    \]The claim follows with $C(\varepsilon)=\max(C_{1}(\varepsilon), c_{2}(\varepsilon))$. Theorem \ref{T4} is proved.\\[-5pt]

    \textsc{Proof of Theorem \ref{T5}.}  Given an integer $a$ and a prime $p$, the congruence $y^{2}\equiv a$ (mod $p$) has at most $2$ solutions, and hence
    \begin{equation}\label{T5.ELL.C.EST.BASIC}
    1 \leq \#E(\mathbb{F}_{p})\leq 1+ 2p
    \end{equation} (the first inequality follows from the fact that at least $\mathcal{O}\in E(\mathbb{F}_{p})$). We have
    \[
    \frac{\#E(\mathbb{F}_{p})}{\varphi(\#E(\mathbb{F}_{p}))}\geq 1
    \]for any prime $p$, and therefore the first inequality in \eqref{T5.Ellip.CURVE.INEQ} is trivial.

    Let us prove the second inequality in \eqref{T5.Ellip.CURVE.INEQ}. We use the following result of David and Wu \cite[Theorem 2.3(i)]{David.Wu}. Suppose that an elliptic curve $E$ does not have complex multiplication. Let $a$ and $t$ be integers with $t\geq 1$. Then we have
    \[
    \#\{p\leq x: \#E(\mathbb{F}_{p})\equiv a\ \text{(mod t)}\}\leq
    C(E) \biggl(\frac{\pi(x)}{\varphi (t)}+ x\, \textup{exp} (-b t^{-2} \sqrt{\ln x})\biggr)
    \]for $\ln x \geq c t^{12} \ln t$. Here $b$ and $c$ are positive absolute constants, and $C(E)>0$ is a constant depending only on $E$.

    We assume that $x \geq c_0 (s)$, where $c_0 (s)>0$ is a constant depending only on $s$; we will choose the constant $c_0 (s)$ later, and it will be large enough. Given a positive integer $t$, we have $t^{12}\ln t \leq t^{13}$. Therefore we obtain
    \begin{equation}\label{David.Wu.Ineq}
     \#\{p\leq x: \#E(\mathbb{F}_{p})\equiv a\ \text{(mod t)}\}\leq
    C(E) \biggl(\frac{\pi(x)}{\varphi (t)}+ x\, \textup{exp} (-b t^{-2} \sqrt{\ln x})\biggr)
    \end{equation}for $1 \leq t \leq (c_1 \ln x)^{1/13}$, where $c_1=1/c$ is a positive absolute constant. We see from \eqref{T5.ELL.C.EST.BASIC} that $\#E(\mathbb{F}_{p})\leq 3p$ for any prime $p$. We put $M=3x$. Hence, $\#E(\mathbb{F}_{p})\leq M$ for any prime $p\leq x$. We have
    \begin{equation}\label{T5.M.RANGE}
    2 \leq (\ln M)^{1/26}\leq (c_1 \ln x)^{1/13},
    \end{equation} if $c_0 (s)$ is chosen large enough.

    Applying Theorem \ref{T1} with $\alpha=1/26$, we have
    \begin{equation}\label{T5.INEQUALITY.FINAL}
    \sum_{p\leq x} \biggl(\frac{\#E(\mathbb{F}_{p})}{\varphi(\#E(\mathbb{F}_{p}))}\biggr)^{s}\leq
    c^{s}\biggl(\pi(x)+ \sum_{q \leq (\ln M)^{1/26}} \frac{\omega(q) (\ln q)^{s}}{q}\biggr),
    \end{equation} where $c$ is a positive absolute constant and
    \[
    \omega (q) = \#\{p\leq x: \#E(\mathbb{F}_{p})\equiv 0\ \text{(mod $q$)}\}.
    \]It follows from \eqref{David.Wu.Ineq} (with $a=0$) and \eqref{T5.M.RANGE} that
    \[
     \omega (q) \leq C(E) \biggl(\frac{\pi(x)}{\varphi (q)}+ x\, \textup{exp} (-b q^{-2} \sqrt{\ln x})\biggr)
    \]for any prime $q \leq (\ln M)^{1/26}$. We obtain
    \begin{align}
    &\sum_{q \leq (\ln M)^{1/26}} \frac{\omega(q) (\ln q)^{s}}{q}\notag\\ 
    &\ \ \ \ \leq
    C(E) \biggl(\pi(x)\sum_{q \leq (\ln M)^{1/26}}\frac{(\ln q)^{s}}{q\varphi (q)}+x\sum_{q \leq (\ln M)^{1/26}}\frac{ (\ln q)^{s}}{q\, \textup{exp} (b q^{-2} \sqrt{\ln x})}\biggr).\label{T5.INEQUALITY.BASIC}
    \end{align}

    Since $\varphi(n) \geq c n/\ln\ln (n+2)$ for any positive integer $n$, where $c$ is a positive absolute constant, we have
    \begin{align}
    \sum_{q \leq (\ln M)^{1/26}}\frac{(\ln q)^{s}}{q\varphi (q)}&\leq \frac{1}{c}
    \sum_{q \leq (\ln M)^{1/26}}\frac{(\ln q)^{s}\ln\ln(q+2)}{q^{2}}\notag\\
    &\leq \frac{1}{c}\sum_{n=1}^{\infty}\frac{(\ln n)^{s}\ln\ln(n+2)}{n^{2}}=c_{1}(s),\label{T5.INEQUALITY.I}
    \end{align} where $c_{1}(s)>0$ is a constant, depending only on $s$.

    We recall that $M=3x$. We have $(\ln (3x))^{1/13}\leq 2 (\ln x)^{1/13}$, if $c_{0}(s)$ is chosen large enough. Hence,
    \[
    b\,\frac{\sqrt{\ln x}}{q^{2}}\geq b\, \frac{(\ln x)^{1/2}}{(\ln (3x))^{1/13}}
    \geq \frac{b}{2}\frac{(\ln x)^{1/2}}{\bigl(\ln x\bigr)^{1/13}}=b_{1} (\ln x)^{11/26}
    \]for any prime $q\leq (\ln (3x))^{1/26}$, where $b_{1}=b/2$ is a positive absolute constant. We obtain
    \[
    x\sum_{q \leq (\ln M)^{1/26}}\frac{ (\ln q)^{s}}{q\, \textup{exp} (b q^{-2} \sqrt{\ln x})}\leq
    x\,\textup{exp}\bigl(-b_{1} (\ln x)^{11/26}\bigr)\sum_{q \leq (\ln M)^{1/26}}\frac{ (\ln q)^{s}}{q}.
    \]

     Putting $k=[(\ln(3x))^{1/26}]$ and applying Lemma \ref{L:primes}, we obtain
    \begin{align*}
    \sum_{q \leq (\ln M)^{1/26}}\frac{ (\ln q)^{s}}{q}&=
    \sum_{q \leq k}\frac{ (\ln q)^{s}}{q} \leq c (\ln k)^{s}\leq
    c \bigl(\ln \bigl((\ln(3x))^{1/26}\bigr)\bigr)^{s}\\
    &=\frac{c}{(26)^{s}}\bigl(\ln\ln (3x)\bigr)^{s},
    \end{align*}where $c$ is a positive absolute constant. We obtain
    \[
    x\sum_{q \leq (\ln M)^{1/26}}\frac{ (\ln q)^{s}}{q\, \textup{exp} (b q^{-2} \sqrt{\ln x})}\leq
    \frac{c}{(26)^{s}}x\,\textup{exp}\bigl(-b_{1} (\ln x)^{11/26}\bigr)\bigl(\ln\ln (3x)\bigr)^{s}.
    \]

    We have $\pi(t) \geq a t/\ln t$ for any real number $t\geq 2$, where $a$ is a positive absolute constant. Let us show that
    \begin{equation}\label{T5.INEQ.FOR.S2}
    \frac{c}{(26)^{s}}x\,\textup{exp}\bigl(-b_{1} (\ln x)^{11/26}\bigr)\bigl(\ln\ln (3x)\bigr)^{s} \leq \frac{a x}{(26)^{s}\ln x}.
    \end{equation} The inequality \eqref{T5.INEQ.FOR.S2} is equivalent to the inequlity
    \[
     c \ln x\bigl(\ln\ln (3x)\bigr)^{s}\leq a\,\textup{exp}\bigl(b_{1} (\ln x)^{11/26}\bigr).
    \]Taking logarithms, we obtain
    \[
    \ln c + \ln\ln x + s \ln\ln\ln (3x)\leq \ln a + b_{1} (\ln x)^{11/26}.
    \]This inequality holds, if $c_{0}(s)$ is chosen large enough. The inequality \eqref{T5.INEQ.FOR.S2} is proved.

    We obtain
    \begin{equation}\label{T5.INEQUALITY.II}
    x\sum_{q \leq (\ln M)^{1/26}}\frac{(\ln q)^{s}}{q\, \textup{exp} (b q^{-2}\sqrt{\ln x})}\leq
    \frac{\pi(x)}{(26)^{s}}.
    \end{equation} Substituting \eqref{T5.INEQUALITY.I} and \eqref{T5.INEQUALITY.II} into \eqref{T5.INEQUALITY.BASIC}, we obtain
    \begin{equation}\label{T5.SUM.OMEGA.ESTIMATE}
    \sum_{q \leq (\ln M)^{1/26}} \frac{\omega(q) (\ln q)^{s}}{q}\leq
    C(E) \biggl(c_{1}(s)+\frac{1}{(26)^{s}}\biggr)\pi(x).
    \end{equation} Substituting \eqref{T5.SUM.OMEGA.ESTIMATE} into \eqref{T5.INEQUALITY.FINAL}, we obtain
    \[
    \sum_{p\leq x} \biggl(\frac{\#E(\mathbb{F}_{p})}{\varphi(\#E(\mathbb{F}_{p}))}\biggr)^{s}\leq
    c^{s}\biggl(1+ C(E) \biggl(c_{1}(s)+\frac{1}{(26)^{s}}\biggr) \biggr)\pi (x) = C_{1}(E,s) \pi(x),
    \]where $C_{1}(E,s)>0$ is a constant depending only on $E$ and $s$. Thus, Theorem \ref{T5} in the case $x\geq c_{0}(s)$ is proved.

    Suppose that $2 \leq x < c_{0}(s)$. We have
    \[
    \#E(\mathbb{F}_{p}) \leq 3 p\leq 3 x\leq 3 c_{0}(s)=c_{2}(s)
    \]for any prime $p\leq x$. Hence,
    \[
    \frac{\#E(\mathbb{F}_{p})}{\varphi(\#E(\mathbb{F}_{p}))}\leq
    \#E(\mathbb{F}_{p})\leq c_{2}(s)
    \]for any prime $p\leq x$. We obtain
    \[
    \sum_{p\leq x} \biggl(\frac{\#E(\mathbb{F}_{p})}{\varphi(\#E(\mathbb{F}_{p}))}\biggr)^{s}\leq
    (c_{2}(s))^{s} \pi (x) = c_{3}(s) \pi (x),
    \]where $c_{3}(s)>0$ is a constant depending only on $s$. The claim follows with $C(E,s)=\max (C_{1}(E,s), c_{3}(s))$. Theorem \ref{T5} is proved.\\[-5pt]

    \textsc{Proof of Theorem \ref{T6}.} Our proof consists of three steps.

    I) We obtain an upper bound for
\[
\sum_{n\leq x} (r(n))^{2}.
\]

We assume that $x\geq x_0$. Since $\textup{ord}_{A}(s)<+\infty$ for any positive integer $s$, we see that $0\leq r(n)<+\infty$ for any positive integer $n$. Since
\[
N_{A}(x)=\sum_{n\leq x}\textup{ord}_{A}(n),
\]we have $N_{A}(x)<+\infty$. By the assumption of the theorem, $N_{A}(x)>0$. Hence, $0< N_{A}(x)<+\infty$. Since $N_{A}(x)>0$, there is a positive integer $n\leq x$ such that $\textup{ord}_{A}(n)>0$. Hence,
\[
0<\rho_{A}(x):= \max_{n\leq x}\textup{ord}_{A}(n)<+\infty.
\]

 It can be shown that
\[
\sum_{n\leq x} (r(n))^{2}=\sum_{\substack{p_1, p_2\in \mathbb{P},\\
j, k\in \mathbb{N}:\\
p_1+a_j\leq x,\\
p_2+a_k\leq x,\\
p_1+a_j=p_2+a_k}}1= \sum_{\substack{\dots\\ p_1=p_2}}1+\sum_{\substack{\dots\\ p_1<p_2}}1+\sum_{\substack{\dots\\ p_1>p_2}}1=
 T_1+ T_2 + T_3.
\]It is easy to see that $T_2=T_3$. Let us estimate $T_1$. We have
\[
T_1= \sum_{\substack{p_1\in \mathbb{P}:\\ p_1\leq x}}\sum_{\substack{j\in \mathbb{N}:\\ a_j\leq x-p_1 }}
\sum_{\substack{k\in \mathbb{N}:\\ a_k=a_j}}1.
\]Since
\[
\sum_{\substack{k\in \mathbb{N}:\\ a_k=a_j}}1=\textup{ord}_{A}(a_j),
\] we obtain
\[
T_1\leq \sum_{\substack{p_1\in \mathbb{P}:\\ p_1\leq x}}\sum_{\substack{j\in \mathbb{N}:\\ a_j\leq x }}
\textup{ord}_{A}(a_j)= \sum_{\substack{j\in \mathbb{N}:\\ a_j\leq x }}
\textup{ord}_{A}(a_j) \sum_{\substack{p_1\in \mathbb{P}:\\ p_1\leq x}} 1=
\pi(x)\sum_{\substack{j\in \mathbb{N}:\\ a_j\leq x }}
\textup{ord}_{A}(a_j).
\]Since $a_j \leq x$, we have $\textup{ord}_{A}(a_j)\leq \rho_{A}(x)$. Hence,
\[
\sum_{\substack{j\in \mathbb{N}:\\ a_j\leq x }}
\textup{ord}_{A}(a_j)\leq \rho_{A}(x)\sum_{\substack{j\in \mathbb{N}:\\ a_j\leq x }}1=
\rho_{A}(x) N_{A}(x).
\]Since, by Chebyshev's theorem, $\pi(x)\leq b x/\ln x,$ where $b>0$ is an absolute constant, we obtain
\[
T_1 \leq b\, \frac{x}{\ln x}\, \rho_{A}(x) N_{A}(x).
\]

Given $a \in \mathbb{N}$, we put
\[
\pi_{2}(x, a)=\#\{p\leq x: \text{$p+a$ is a prime}\}.
\]We use the following result of Schnirelmann \cite{Schnirelmann}. Let $a$ be a positive integer, $x$ be a real number with $x\geq 4$. Then
\[
\pi_{2}(x, a)\leq c\,\frac{x}{(\ln x)^2}\,\frac{a}{\varphi(a)},
\]where $c>0$ is an absolute constant. Let us estimate $T_2$. We have
\begin{align*}
T_2&= \sum_{\substack{p_1, p_2\in \mathbb{P},\\
j, k\in \mathbb{N}:\\
p_1+a_j\leq x,\\
p_2+a_k\leq x,\\
p_1+a_j=p_2+a_k,\\
p_1<p_2}}1
\leq \sum_{\substack{j, k\in \mathbb{N}:\\ a_k < a_j\leq x}}
\sum_{\substack{p\in \mathbb{P}:\\p\leq x\\ \text{$p+a_j-a_k$ is a prime}}}1= \sum_{\substack{j, k\in \mathbb{N}:\\ a_k < a_j\leq x}}\pi_{2}(x, a_j - a_k)\\
&\leq
c\,\frac{x}{(\ln x)^2}\sum_{\substack{j, k\in \mathbb{N}:\\ a_k < a_j\leq x}}\frac{a_j - a_k}{\varphi(a_j - a_k)}= c\,\frac{x}{(\ln x)^2}\sum_{\substack{k\in \mathbb{N}:\\ a_k < x}}\sum_{\substack{j\in \mathbb{N}:\\ a_k < a_j\leq x}}\frac{a_j - a_k}{\varphi(a_j - a_k)}.
\end{align*} Fix a positive integer $k$ with $a_k < x$. Let $j$ be a positive integer such that $a_k < a_j\leq x$.
Then $0< a_j - a_k\leq x$. Applying Theorem \ref{T1} with $s=1$, $M=x$, and $\alpha$, given by the assumption of the theorem, we obtain
\[
\sum_{\substack{j\in \mathbb{N}:\\ a_k < a_j\leq x}}\frac{a_j - a_k}{\varphi(a_j - a_k)}
\leq C(\alpha)\biggl(l+\sum_{p\leq (\ln x)^{\alpha}}
\frac{\omega(p)
\ln p}{p}\biggr),
\]where
\[
l=\#\{j\in \mathbb{N}: a_k<a_j\leq x\}\leq \#\{j\in \mathbb{N}: a_j\leq x\}=N_{A}(x)
\]and
\[
\omega(p)=\#\{j\in \mathbb{N}: \text{$a_k<a_j\leq x$ and $a_j\equiv a_k$ (mod $p$)}\}.
\] Hence,
\begin{align*}
&\sum_{\substack{j\in \mathbb{N}:\\ a_k < a_j\leq x}}\frac{a_j - a_k}{\varphi(a_j - a_k)}\\
&\quad \ \,\,\leq C(\alpha)\biggl(N_{A}(x)+\sum_{p\leq (\ln x)^{\alpha}}\frac{\#\{j\in \mathbb{N}: \text{$a_k<a_j\leq x$ and $a_j\equiv a_k$ (mod $p$)}\}
\ln p}{p}\biggr).
\end{align*} We obtain
\begin{align*}
T_2 &\leq c\,\frac{x}{(\ln x)^2}\sum_{\substack{k\in \mathbb{N}:\\ a_k < x}}C(\alpha)\bigg(N_{A}(x)\\
&\quad\ +\sum_{p\leq (\ln x)^{\alpha}}\frac{\#\{j\in \mathbb{N}: \text{$a_k<a_j\leq x$ and $a_j\equiv a_k$ (mod $p$)}\}
\ln p}{p}\bigg)\\
&\leq c C(\alpha)\,\frac{x}{(\ln x)^2}(N_{A}(x))^{2}+c C(\alpha)\,\frac{x}{(\ln x)^2}\\
&\quad\ \times\sum_{\substack{k\in \mathbb{N}:\\ a_k < x}}\sum_{p\leq (\ln x)^{\alpha}}\frac{\#\{j\in \mathbb{N}: \text{$a_k<a_j\leq x$ and $a_j\equiv a_k$ (mod $p$)}\}\ln p}{p}.
\end{align*}By the assumption of the theorem,
\[
\sum_{\substack{k\in \mathbb{N}:\\ a_k < x}}\sum_{p\leq (\ln x)^{\alpha}}\frac{\#\{j\in \mathbb{N}: \text{$a_k<a_j\leq x$ and $a_j\equiv a_k$ (mod $p$)}\}\ln p}{p}\leq \gamma_{2}(N_{A}(x))^{2}.
\]We obtain
\[
T_2\leq c_{0}(\gamma_2, \alpha)\,\frac{x}{(\ln x)^2}(N_{A}(x))^{2},
\]where $c_{0}(\gamma_2,\alpha)>0$ is a constant, depending only on $\gamma_2$ and $\alpha$.

We obtain
\begin{align*}
\sum_{n\leq x} (r(n))^{2}&= T_1+T_2+T_3=T_1+2T_2\\
&\leq b\, \frac{x}{\ln x}\, \rho_{A}(x) N_{A}(x)  + 2c_{0}(\gamma_2, \alpha)\,\frac{x}{(\ln x)^2}(N_{A}(x))^{2}\\
&\leq c(\gamma_2, \alpha)\,\frac{x}{(\ln x)^2}\, N_{A}(x)\bigl(\rho_{A}(x)\ln x+ N_{A}(x)\bigr),
\end{align*}where $c(\gamma_2, \alpha)=b+2c_{0}(\gamma_2, \alpha)>0$ is a constant, depending only on $\gamma_2$ and $\alpha$.

II) We obtain a lower bound for
\[
\sum_{\substack{n\leq x:\\ r(n)\geq b_1 \gamma_1 (N_{A}(x)/\ln x)}} r(n),
\]where $b_1>0$ is some absolute constant, which will be defined later.

By the assumption of the theorem, $N_{A}(x/2)\geq \gamma_1 N_{A}(x)>0$. Also we have $N_{A}(x/2)<+\infty$. Hence, $0< N_{A}(x/2)<+\infty$. Since $x\geq x_0\geq 10$, we have
\[
\pi(x/2)\geq b\,\frac{x/2}{\ln x/2}\geq\frac{b}{2}\,\frac{x}{\ln x}=
b_0\, \frac{x}{\ln x},
\]where $b_0>0$ is an absolute constant.

We obtain
\[
\sum_{n\leq x} r(n)\geq \pi \bigg(\frac{x}{2}\bigg)N_{A}\bigg(\frac{x}{2}\bigg)\geq b_{0}\gamma_{1}\,\frac{x}{\ln x}\,N_{A}(x).
\]Also,
\begin{align*}
\sum_{\substack{n\leq x:\\ r(n)< (b_{0}\gamma_{1}N_{A}(x))/ (2\ln x)}} r(n) &< \frac{b_{0}\gamma_{1}}{2}\,\frac{N_{A}(x)}{\ln x}
\sum_{\substack{n\leq x:\\ r(n)< (b_{0}\gamma_{1}N_{A}(x))/ (2\ln x)}} 1\\
&\leq \frac{b_{0}\gamma_{1}}{2}\,\frac{N_{A}(x)}{\ln x} \sum_{n\leq x} 1 \leq \frac{b_{0}\gamma_{1}}{2}\,\frac{x}{\ln x}\,N_{A}(x).
\end{align*} Hence,
\begin{align*}
\sum_{\substack{n\leq x:\\ r(n)\geq (b_{0}\gamma_{1}N_{A}(x))/ (2\ln x)}} r(n) &=
\sum_{n\leq x} r(n) - \sum_{\substack{n\leq x:\\ r(n)< (b_{0}\gamma_{1}N_{A}(x))/ (2\ln x)}} r(n)\\
&\geq \frac{b_{0}\gamma_{1}}{2}\,\frac{x}{\ln x}\,N_{A}(x).
\end{align*} We put $b_1=b_{0}/2$. Then $b_1>0$ is an absolute constant, and
\[
\sum_{\substack{n\leq x:\\ r(n)\geq b_1\gamma_1 (N_{A}(x)/\ln x)}} r(n)\geq b_1 \gamma_1 \,\frac{x}{\ln x}\, N_{A}(x).
\]

III) Applying the Cauchy-Schwarz inequality, we have
\begin{align*}
&\bigg(b_1 \gamma_1 \,\frac{x}{\ln x}\, N_{A}(x)\bigg)^2 \leq \Bigg(\sum_{\substack{n\leq x:\\ r(n)\geq b_1\gamma_1 (N_{A}(x)/\ln x)}} r(n)\Bigg)^{2}\\
&\leq \Bigg(\sum_{\substack{n\leq x:\\ r(n)\geq b_1\gamma_1 (N_{A}(x)/\ln x)}} 1\Bigg)\cdot
\Bigg(\sum_{\substack{n\leq x:\\ r(n)\geq b_1\gamma_1 (N_{A}(x)/\ln x)}} (r(n))^{2} \Bigg)\\
&\leq \Bigg(\sum_{\substack{n\leq x:\\ r(n)\geq b_1\gamma_1 (N_{A}(x)/\ln x)}} 1\Bigg)\cdot
\Bigg(\sum_{n\leq x} (r(n))^{2} \Bigg)\\
&\leq \Bigg(\sum_{\substack{n\leq x:\\ r(n)\geq b_1\gamma_1 (N_{A}(x)/\ln x)}} 1\Bigg)\cdot  c(\gamma_2, \alpha)\,\frac{x}{(\ln x)^2}\, N_{A}(x)\bigl(\rho_{A}(x)\ln x+ N_{A}(x)\bigr).
\end{align*}We obtain
\[
 \sum_{\substack{n\leq x:\\ r(n)\geq b_1\gamma_1 (N_{A}(x)/\ln x)}} 1\geq
  \frac{(b_1 \gamma_1)^2}{c(\gamma_2, \alpha)}\,x\,\frac{N_{A}(x)}{N_{A}(x)+\rho_{A}(x)\ln x}.
\]We put
\[
c_1= b_1\gamma_1,\qquad c_2= \frac{(b_1 \gamma_1)^2}{c(\gamma_2, \alpha)}.
\]Then $c_1$ and $c_2$ are positive constants, depending only on $\gamma_1$ and $\gamma_1$, $\gamma_2$, $\alpha$ respectively, and
\[
\#\left\{1\leq n \leq x: r(n)\geq c_1\, \frac{N_{A}(x)}{\ln x}\right\}\geq c_2 x\,\frac{N_{A}(x)}{N_{A}(x)+\rho_{A}(x)\ln x}.
\] Theorem \ref{T6} is proved.\\[-5pt]

\textsc{Proof of Theorem \ref{T7}.} We are going to apply Theorem \ref{T6}. We put
\[
\Omega = \{n\in \mathbb{N}: R(n)>0\},\qquad A= \{R(n): n\in \Omega\}.
\] Given a positive integer $m$, we have
\[
\textup{ord}_{A}(m)=\#\{n\in \Omega: R(n)=m\}\leq k.
\] In particular, we see that $\textup{ord}_{A}(m) <+\infty$ for any positive integer $m$, and
\[
\rho_{A}(t)= \max_{n\leq t} \textup{ord}_{A}(n) \leq k
\]for any real number $t\geq 1$.

There is a positive integer $N_{0}$, depending only on the polynomial $R$, such that
\[
-\frac{a_{k}}{2} t^{k} \leq a_{k-1}t^{k-1}+\ldots+a_{0}\leq \frac{a_{k}}{2} t^{k}
\]and
\[
R'(t)=k a_k t^{k-1}+ (k-1)a_{k-1}t^{k-2}+\ldots + a_1> 0
\]
 for any real number $t \geq N_{0}$. Hence, $(a_{k}/2) n^{k} \leq R(n) \leq 2 a_{k} n^{k}$ and $R(n)< R(n+1)$ for any integer $n \geq N_{0}$. We put
\[
\widetilde{M}_{2}=\max_{\substack{n\in \mathbb{N}:\\ n\leq N_{0}}} R(n).
\]We see that $\widetilde{M}_{2}$ is a positive constant, depending only on $R$, and
\[
R(n)\leq \widetilde{M}_{2}\leq \widetilde{M}_{2} n^{k}
\] for any integer $n$ such that $n\in \Omega$ and $n\leq N_{0}$. Therefore
\[
R(n) \leq \max(\widetilde{M}_{2}, 2a_{k})n^{k} = M_{2} n^{k}
\] for any $n\in \Omega$.

We put
\[
\widetilde{M}_{1}= \frac{1}{(N_{0})^{k}}.
\]Hence,
\[
\widetilde{M}_{1} n^{k} \leq \widetilde{M}_{1} (N_{0})^{k}=1 \leq R(n)
\] for any integer $n$ such that $n\in \Omega$ and $n \leq N_{0}$. We obtain
\[
R(n) \geq \min (\widetilde{M}_{1}, a_{k}/2) n^{k}= M_{1} n^{k}
\]for any $n\in \Omega$. Thus,
\begin{equation}\label{T7.M1.R.M2}
M_{1} n^{k} \leq R(n) \leq M_{2} n^{k}
\end{equation}for any $n\in \Omega$, where $M_{1}$ and $M_{2}$ are positive constants, depending only on $R$.

We have
\begin{equation}\label{T8.OMEGA}
\Omega= \{n_1,\ldots, n_{T}, N_{0}, N_{0}+1,\ldots\},
\end{equation} where $n_{1}, \ldots, n_{T}$ are positive integers with $n_{1}<\ldots < n_{T}< N_{0}$. We may assume that $T>0$. We see that $T$ is a positive constant, depending only on $R$.

We assume that $x  \geq x_{0}$, where $x_{0}$ is a positive constant, depending only on $R$; we will choose the constant $x_{0}$ later, and it will be large enough. Let $x_{0} \geq M_{2} (N_{0})^{k}$. Applying \eqref{T7.M1.R.M2}, we have
\begin{align*}
N_{A}(x)&=\#\{n\in \Omega: R(n)\leq x\}\leq
\#\{n\in \Omega: M_{1} n^{k}\leq x\}\\
&\leq \#\biggl\{n\in \mathbb{N}: n\leq \bigg(\frac{x}{M_{1}}\bigg)^{1/k}\biggr\}=
\biggl[\bigg(\frac{x}{M_{1}}\bigg)^{1/k}\biggr]\leq \bigg(\frac{x}{M_{1}}\bigg)^{1/k}.
\end{align*}

Let $x_1 = M_{2} (2N_{0}+1)^{k}$. Then $x_{1}$ is a positive constant depending only on $R$. Let $t$ be a real number such that $t \geq x_{1}$. Applying \eqref{T7.M1.R.M2} and \eqref{T8.OMEGA}, we have
\begin{align}
N_{A}(t)&= \#\{n\in \Omega: R(n)\leq t\}\geq
\#\{n\in \Omega: M_{2}n^{k}\leq t\}\notag\\
&= \#\bigg\{n\in \Omega: n\leq \bigg(\frac{t}{M_{2}}\bigg)^{1/k}\bigg\}\geq \biggl[\bigg(\frac{t}{M_{2}}\bigg)^{1/k}\biggr] - N_{0}+1\notag \\
&\geq \biggl(\frac{t}{M_{2}}\biggr)^{1/k} - N_{0}\geq \frac{1}{2} \biggl(\frac{t}{M_{2}}\biggr)^{1/k} .\label{T7.N.A.t.LOW.BOUND}
\end{align}

Let $x_{0} \geq 2x_1$. Then
\begin{equation}\label{T7.N.A.x.BASIC.EST}
\frac{1}{2} \biggl(\frac{x}{M_{2}}\biggr)^{1/k} \leq N_{A}(x) \leq \biggl(\frac{x}{M_{1}}\biggr)^{1/k}.
\end{equation}In particular, we see that $N_{A}(x)>0$. Since $x/2 \geq x_{1}$, from \eqref{T7.N.A.t.LOW.BOUND} we have
\[
N_{A}(x/2) \geq \frac{1}{2} \biggl(\frac{x/2}{M_{2}}\biggr)^{1/k}=
\frac{1}{2} \biggl(\frac{M_{1}}{2M_{2}}\biggr)^{1/k} \biggl(\frac{x}{M_{1}}\biggr)^{1/k}\geq \frac{1}{2} \biggl(\frac{M_{1}}{2M_{2}}\biggr)^{1/k} N_{A}(x)= \gamma_{1} N_{A}(x),
\]where $\gamma_{1}$ is a positive constant, depending only on $R$. We see that \eqref{T6.INEQ1} and \eqref{T6.INEQ2} hold.

We may assume that $x_{0}\geq R(N_{0})$. We have
\[
\rho_{A}(x) = \max_{n\leq x} \textup{ord}_{A}(n)\geq \textup{ord}_{A}\bigl(R(N_{0})\bigr) \geq 1.
\]We obtain
\begin{equation}\label{T7.rho}
1 \leq \rho_{A}(x) \leq k .
\end{equation}

Let us consider the sum
\[
S= \sum_{\substack{j\in \Omega:\\ R(j)< x}} \sum_{p\leq \ln x} \frac{\lambda (j, p)\ln p}{p},
\]where
\[
\lambda (j, p)= \#\{n\in \Omega: R(j)< R(n)\leq x \text{ and $R(n) \equiv R(j)$ (mod $p$)}\}.
\] Fix $j$ from the range of summation of $S$. We have
\begin{equation}\label{T7.lambda}
\sum_{p\leq \ln x} \frac{\lambda (j, p)\ln p}{p}= \sum_{\substack{p\leq \ln x\\ p|a_k}} \frac{\lambda (j, p)\ln p}{p}+
\sum_{\substack{p\leq \ln x\\ (p,a_k)=1}} \frac{\lambda (j, p)\ln p}{p}= S_{1}+ S_{2}.
\end{equation}

Let $p$ be in the range of summation of $S_{1}$. We have $\lambda (j,p) \leq N_{A}(x)$. We obtain
\begin{align}\label{T7.S1}
S_{1} \leq N_{A}(x)\sum_{\substack{p\leq \ln x\\ p|a_k}} \frac{\ln p}{p}\leq \biggl(\sum_{p|a_k} \frac{\ln p}{p} + 1\biggr)N_{A}(x)= c_{1} N_{A}(x),
\end{align}where $c_1$ is a positive constant, depending only on $R$.

By Bertrand's postulate (see, for example, \cite[Theorem 3.1.9]{Murty}), there is a positive integer $n_0$ such that there is a prime between $n$ and $2 n$ for any integer $n \geq n_0$. We see that there is a prime $p$ between $n_{0}+a_k$ and $2(n_{0}+a_k)$. We have $p> a_k$, and hence $(p, a_k)=1$. We have $\ln x > 2(n_{0}+a_k)$, if $x_0$ is chosen large enough. We see that the set $\{p: p\leq \ln x \text{ and } (p,a_k)=1\}\neq\emptyset$.

We may assume that (see \eqref{T8.OMEGA})
\[
x_{0}\geq \max \bigl(R(n_1),\ldots, R(n_{T}), R(N_{0}),\ldots, R(N_{0}+10)\bigr).
\]We define
\[
\Omega(x)= \{n\in \Omega: R(n)\leq x\}.
\]Since $R(n)< R(n+1)$ for any integer $n \geq N_{0}$, we have
\begin{equation}\label{T7.OMEGA.x.S}
\Omega(x)=\{n_{1},\ldots, n_{T}, N_{0}, N_{0}+1,\ldots, N_{0}+r\}.
\end{equation}

Let $p$ be in the range of summation of $S_{2}$. We define
\[
U=\bigl\{b\in \{0,\ldots, p-1\}: R(b)\equiv R(j) \text{ (mod $p$)}\bigr\}.
\] Trivially, $\#U \leq p$. Since $(p,a_k)=1$, we have $\#U \leq k$. We obtain $\#U \leq \min (p,k)$. We observe that if $b\in \{0,\ldots, p-1\}$ such that $b\equiv j$ (mod $p$), then $b\in U$. Therefore
\[
1\leq \# U\leq \min (p,k).
\]

Given $b\in U$, we define
\[
\Lambda (b)= \{ t\in \mathbb{Z}: b+pt\in \Omega(x)\}.
\] We see from \eqref{T7.OMEGA.x.S} that
\begin{align*}
\# \Lambda (b)&\leq T + \#\{t\in \mathbb{Z}: N_{0}\leq b+ pt\leq N_{0}+r\}\\
&= T + \biggl[\frac{N_{0}+r - b}{p}\biggr] - \biggl\lceil\frac{N_{0}-b}{p}\biggr\rceil + 1\leq \frac{r}{p}+ T+1.
\end{align*} We have
\[
r\leq \#\Omega(x) = N_{A}(x).
\]We have (see \eqref{T7.N.A.x.BASIC.EST})
\[
N_{A}(x) \geq \frac{1}{2} \biggl(\frac{x}{M_{2}}\biggr)^{1/k} \geq \ln x,
\]if $x_0$ is chosen large enough. Since $p\leq \ln x$, we obtain $p \leq N_{A}(x)$. Therefore
\[
\# \Lambda (b)\leq \frac{N_{A}(x)}{p}+ T+1\leq \frac{N_{A}(x)}{p}+ (T+1)\frac{N_{A}(x)}{p}= (T+2)\frac{N_{A}(x)}{p}= c_2 \frac{N_{A}(x)}{p},
\]where $c_2$ is a positive constant, depending only on $R$.

We have
\[
\lambda (j,p) \leq \sum_{b\in U} \# \Lambda (b)\leq c_2\,\frac{N_{A}(x)}{p}\,\#U\leq
c_2 k\,\frac{N_{A}(x)}{p} = c_3\,\frac{N_{A}(x)}{p},
\]where $c_3$ is a positive constant, depending only on $R$. We obtain
\begin{equation}
S_2=\sum_{\substack{p\leq \ln x\\ (p,a_k)=1}} \frac{\lambda (j, p)\ln p}{p}\leq
c_3 N_{A}(x) \sum_{\substack{p\leq \ln x\\ (p,a_k)=1}} \frac{\ln p}{p^{2}}\leq c_3 N_{A}(x) \sum_{p} \frac{\ln p}{p^{2}} = c_4 N_{A}(x),\label{T7.S2}
\end{equation} where $c_{4}$ is a positive constant, depending only on $R$.

Substituting \eqref{T7.S1} and \eqref{T7.S2} into \eqref{T7.lambda}, we obtain
\[
\sum_{p\leq \ln x} \frac{\lambda (j, p)\ln p}{p}\leq (c_1 + c_4) N_{A}(x) = \gamma_{2} N_{A}(x),
\]where $\gamma_{2}$ is a positive constant, depending only on $R$. We have
\[
\sum_{\substack{j\in \Omega:\\ R(j)< x}} \sum_{p\leq \ln x} \frac{\lambda (j, p)\ln p}{p}\leq
\gamma_{2} N_{A}(x)\sum_{\substack{j\in \Omega:\\ R(j)< x}} 1 \leq \gamma_{2} \bigl(N_{A}(x)\bigr)^{2}.
\] We see that \eqref{T6.INEQ3} with $\alpha = 1$ holds.

By Theorem \ref{T6}, there are positive constants $c_{1}=c_{1}(\gamma_1)$ and $c_{2}=c_{2}(\gamma_1, \gamma_2, \alpha)$ such that
\[
\#\biggl\{1\leq n \leq x: \widetilde{r}(n)\geq c_{1}\frac{N_{A}(x)}{\ln x}\biggr\}\geq c_{2}x\,\frac{N_{A}(x)}{N_{A}(x)+\rho_{A}(x)\ln x},
\]where
\[
\widetilde{r}(n) = \#\{(p,j)\in \mathbb{P}\times \Omega: p+ R(j)=n\}.
\]Since $\gamma_{1}$ and $\gamma_2$ are positive constants depending only on $R$, $\alpha = 1$, we see that $c_1$ and $c_2$ are positive constants depending only on $R$. Applying \eqref{T7.N.A.x.BASIC.EST}, we obtain
\[
\#\biggl\{1\leq n \leq x: \widetilde{r}(n)\geq c_{1}\frac{N_{A}(x)}{\ln x}\biggr\}\leq \#\biggl\{1\leq n \leq x:  \widetilde{r}(n)\geq \frac{c_{1}}{2 (M_{2})^{1/k}}\frac{x^{1/k}}{\ln x}\biggr\}.
\]We put
\[
r(n) = \#\{(p,j)\in \mathbb{P}\times\mathbb{N}: p+ R(j)=n\}.
\]It is clear that $\widetilde{r}(n) \leq r(n)$ for any positive integer $n$. Hence,
\[
 \#\biggl\{1\leq n \leq x: \widetilde{r}(n)\geq \frac{c_{1}}{2 (M_{2})^{1/k}}\frac{x^{1/k}}{\ln x}\biggr\}\leq
  \#\biggl\{1\leq n \leq x: r(n)\geq \frac{c_{1}}{2 (M_{2})^{1/k}}\frac{x^{1/k}}{\ln x}\biggr\}.
\] We have (see \eqref{T7.rho})
\[
0< \rho_{A}(x)\ln x\leq k \ln x \leq \frac{x^{1/k}}{2 (M_{2})^{1/k}},
\]if $x_0$ is chosen large enough. We obtain $0< \rho_{A}(x)\ln x\leq N_{A}(x)$, and hence
\[
\frac{N_{A}(x)}{N_{A}(x)+\rho_{A}(x)\ln x}\geq \frac{N_{A}(x)}{N_{A}(x)+ N_{A}(x)}=\frac{1}{2}.
\]We obtain
\[
\#\biggl\{1\leq n \leq x: r(n)\geq \frac{c_{1}}{2 (M_{2})^{1/k}}\frac{x^{1/k}}{\ln x}\biggr\}\geq \frac{c_{2}}{2}\,x.
\] We see that $c_1 / (2 (M_{2})^{1/k})$ and $c_2/2$ are positive constants, depending only on $R$. Let us denote $c_1 / (2 (M_{2})^{1/k})$ by $c_1$ and $c_2/2$ by $c_2$. Theorem \ref{T7} is proved.\\[-5pt]

\textsc{Proof of Corollary \ref{Cor2}.} We put $R(n)=n^{k}$. By Theorem \ref{T7}, there are positive constants $c_{1}(k)$, $c_{2}(k)$, and $x_{0}(k)$, depending only on $k$, such that
\[
\#\biggl\{1\leq n \leq x: r(n)\geq c_{1}(k)\frac{x^{1/k}}{\ln x}\biggr\}\geq c_{2}(k) x
\]for any real number $x \geq x_{0}(k)$, where
\[
r(n)=\#\{(p, j)\in \mathbb{P}\times\mathbb{N}: p+j^{k}=n\}.
\] We put
\[
\max_{4 \leq x \leq x_{0}(k)}\frac{x^{1/k}}{\ln x}=\alpha (k).
\]We see that $\alpha(k)$ is a positive constant, depending only on $k$. Since $3=2+1=2+R(1)$ and $2\in \mathbb{P}$, we have $r(3) \geq 1$. We obtain
\[
\#\biggl\{1\leq n \leq x: r(n)\geq \frac{1}{\alpha (k)}\frac{x^{1/k}}{\ln x}\biggr\}\geq 1 \geq \frac{1}{x_{0}(k)} x
\]for any real number $3 \leq x \leq x_{0}(k)$. We put
\[
b_{1}(k) = \min \biggl(c_{1}(k), \frac{1}{\alpha(k)}\biggr),\qquad b_{2}(k)=\min \biggl(c_{2}(k), \frac{1}{x_{0}(k)}\biggr).
\] We see that $b_{1}(k)$ and $b_{2}(k)$ are positive constants, depending only on $k$. We have
\[
\#\biggl\{1\leq n \leq x: r(n)\geq b_{1}(k)\frac{x^{1/k}}{\ln x}\biggr\}\geq b_{2}(k) x
\]for any real number $x\geq 3$. Let us denote $b_{1}(k)$ by $c_{1}(k)$ and $b_{2}(k)$ by $c_{2}(k)$. Corollary \ref{Cor2} is proved.\\[-5pt]

\textsc{Proof of Theorem \ref{T8}.} By Hasse's theorem \cite{Hasse}, we have $\bigl|\#E(\mathbb{F}_{p})- (p+1)\bigr|< 2\sqrt{p}$ for any prime $p> 3$.  We see from \eqref{T5.ELL.C.EST.BASIC} that $|\#E(\mathbb{F}_{p}) - (p+1)|<2\sqrt{p}$ for $p\in \{2, 3\}$. We obtain
\[
|\#E(\mathbb{F}_{p}) - (p+1)|<2\sqrt{p}
\] for any prime $p$. Thus,
\begin{equation}\label{T8.NUMBER.ELL.C.SQRT}
(\sqrt{p}-1)^{2}< \#E(\mathbb{F}_{p})< (\sqrt{p}+1)^{2}
\end{equation} for any prime $p$.

We put $A=\bigl\{ \#E(\mathbb{F}_{p}):\ p \geq 2\bigr\}$. We are going to apply Theorem \ref{T6}. Suppose that $p\geq 5$. Let $q$ be a prime with $q> p + 4\sqrt{p}+ 4$. Then $(\sqrt{q}-1)^{2}> (\sqrt{p}+1)^{2}$, and from \eqref{T8.NUMBER.ELL.C.SQRT} we obtain $\#E(\mathbb{F}_{q})> \#E(\mathbb{F}_{p})$. Let $q$ be a prime with $q< p - 4\sqrt{p}+ 4$. Then $(\sqrt{q}+1)^{2}< (\sqrt{p}-1)^{2}$, and it follows from \eqref{T8.NUMBER.ELL.C.SQRT} that $\#E(\mathbb{F}_{q})< \#E(\mathbb{F}_{p})$. Hence,
\begin{align*}
\textup{ord}_{A}\bigl(\#E(\mathbb{F}_{p})\bigr)&=\#\{q: \#E(\mathbb{F}_{q})=\#E(\mathbb{F}_{p})\}\\
&\leq \#\{q: p - 4\sqrt{p}+ 4 \leq q \leq p + 4\sqrt{p}+ 4 \}\\
&\leq \#\{n\in \mathbb{N}: p - 4\sqrt{p}+ 4 \leq n \leq p + 4\sqrt{p}+ 4 \}\\
&= \bigl[p + 4\sqrt{p}+ 4\bigr] - \bigl\lceil p - 4\sqrt{p}+ 4\bigr\rceil + 1\leq 8\sqrt{p}+1 < 9\sqrt{p}.
\end{align*}Suppose that $p<5$, i.e. $p\in \{2,3\}$. From \eqref{T5.ELL.C.EST.BASIC} we obtain $\#E(\mathbb{F}_{p})\leq 7$. Let $q$ be a prime with $q > 14$. Then $(\sqrt{q}-1)^{2}>7$, and from \eqref{T8.NUMBER.ELL.C.SQRT} we obtain $\#E(\mathbb{F}_{q})> \#E(\mathbb{F}_{p})$. Hence,
\[
\textup{ord}_{A}\bigl(\#E(\mathbb{F}_{p})\bigr)\leq \#\{q: q\leq 14\} = 6 <
6\sqrt{p}.
\] We obtain
\begin{equation}\label{T8.ord.ESTIMATE}
\textup{ord}_{A}\bigl(\#E(\mathbb{F}_{p})\bigr) \leq 9\sqrt{p}
\end{equation} for any prime $p$. In particular, we see that $\textup{ord}_{A} (n) <+\infty$ for any positive integer $n$.

We assume that $x\geq x_0$, where $x_0$ is a positive absolute constant; we will choose the constant $x_0$ later, and it will be large enough. We may assume that $x_{0} \geq 100$, and hence $x\geq 100$. Let $p$ be a prime such that $\#E(\mathbb{F}_{p})\leq x$. It follows from \eqref{T8.NUMBER.ELL.C.SQRT} that $(\sqrt{p}-1)^{2}< x$ or, that is equivalent, $\sqrt{p}-1< \sqrt{x}$. We obtain
\begin{equation}\label{T8.p.EST}
p< x+2\sqrt{x}+1< 2x.
\end{equation}

We put $n_0 = \#E(\mathbb{F}_2)$. We see from \eqref{T5.ELL.C.EST.BASIC} that $1\leq n_0 \leq 5< x$. We obtain
\[
\rho_{A}(x)=\max_{n\leq x}\textup{ord}_{A}(n) \geq \textup{ord}_{A}(n_0) \geq 1.
\] Let $n_1$ be a positive integer such that $n_1 \leq x$ and $\rho_{A}(x)= \textup{ord}_{A}(n_1)$. Since $\rho_{A}(x) \geq 1$, we obtain $\textup{ord}_{A}(n_1) \geq 1$, and hence there is a prime $p$ such that $\#E(\mathbb{F}_{p})=n_1$. Since $n_1\leq x$, we have $\#E(\mathbb{F}_{p}) \leq x$. From \eqref{T8.p.EST} we obtain $p\leq 2x$. Applying \eqref{T8.ord.ESTIMATE}, we have
\[
\rho_{A}(x)= \textup{ord}_{A}(n_1)= \textup{ord}_{A}\bigl(\#E(\mathbb{F}_{p})\bigr)\leq 9\sqrt{p}
\leq 9\sqrt{2x}< 13 \sqrt{x}.
\]Thus,
\begin{equation}\label{T8.rho.EST}
1 \leq \rho_{A}(x) \leq 13\sqrt{x}.
\end{equation}

Applying \eqref{T8.p.EST}, we have
\begin{align*}
N_{A}(x)&=\#\{p: \#E(\mathbb{F}_{p})\leq x\}\leq \#\{p: p\leq 2x\}= \pi (2x)\\
&\leq a\,\frac{2x}{\ln (2x)}\leq 2a\,\frac{x}{\ln x}= a_{2}\,\frac{x}{\ln x},
\end{align*}where $a_2$ is a positive absolute constant.

Let $t$ be a real number with $t\geq 20$. It is easy to see that
\[
\frac{t}{2}\leq t - 2\sqrt{t}+1.
\]Hence, if $p$ is a prime such that $p\leq t/2$, then $p\leq t-2\sqrt{t}+1$ or, that is equivalent, $(\sqrt{p}+1)^{2}\leq t$. From \eqref{T8.NUMBER.ELL.C.SQRT} we obtain $\#E(\mathbb{F}_{p})\leq t$. Therefore
\begin{align}
N_{A}(t)&=\#\{p: \#E(\mathbb{F}_{p})\leq t\}\geq \#\{p: p \leq t/2\}= \pi \biggl(\frac{t}{2}\biggr)\notag \\
&\geq b\,\frac{t/2}{\ln (t/2)}\geq \frac{b}{2}\,\frac{t}{\ln t}= a_{1}\,\frac{t}{\ln t},\label{T8.N.t.est}
\end{align}where $a_1$ is a positive absolute constant. Thus,
\begin{equation}\label{T8.N.A.ESTIM}
a_{1} \frac{x}{\ln x} \leq N_{A}(x) \leq a_{2} \frac{x}{\ln x},
\end{equation}where $a_1$, $a_2$ are positive absolute constants. Since $x/2 >20$, from \eqref{T8.N.t.est} we obtain
\[
N_{A}\bigg(\frac{x}{2}\bigg) \geq a_{1}\,\frac{x/2}{\ln (x/2)}\geq \frac{a_1}{2}\,\frac{x}{\ln x}=
\frac{a_1}{2a_2}\,a_{2}\,\frac{x}{\ln x}\geq \frac{a_1}{2a_2}\,N_{A}(x)=\gamma_1 N_{A}(x),
\]where $\gamma_1 = a_{1}/ (2a_{2})$ is a positive absolute constant. We see that \eqref{T6.INEQ1} and \eqref{T6.INEQ2} hold.

 We have
\[
10 \leq (\ln x)^{1/14} \leq (c_{1}\ln x)^{1/13},
\] if $x_0$ is chosen large enough (here $c_1$ is a positive absolute constant from \eqref{David.Wu.Ineq}). Let us consider the sum
\[
S=\sum_{\substack{q\in \mathbb{P}:\\ \#E(\mathbb{F}_{q})< x}}\sum_{\substack{t\in \mathbb{P}:\\ t\leq (\ln x)^{1/14}}} \frac{\lambda (q,t) \ln t}{t},
\]where
\[
\lambda (q, t):= \#\{p\in \mathbb{P}: \#E(\mathbb{F}_{q})< \#E(\mathbb{F}_{p})\leq x \text{ and $\#E(\mathbb{F}_{p})\equiv \#E(\mathbb{F}_{q})$ (mod $t$)}\}.
\]Let $q$ and $t$ be in the range of summation of $S$. From \eqref{T8.p.EST} and \eqref{David.Wu.Ineq} we obtain
\begin{align*}
\lambda (q, t) &\leq \#\{p\leq 2x: \#E(\mathbb{F}_{p}) \equiv \#E(\mathbb{F}_{q}) \text{ (mod $t$)}\}\\
&\leq C(E) \biggl(\frac{\pi(2x)}{\varphi (t)}+ 2x\, \textup{exp} (-b t^{-2} \sqrt{\ln (2x)})\biggr).
\end{align*}Hence,
\begin{align}
&\sum_{\substack{t\in \mathbb{P}:\\ t\leq (\ln x)^{1/14}}} \frac{\lambda (q,t) \ln t}{t}\notag\\
&\quad\ \ \leq
C(E) \Biggl(\pi(2x)\sum_{\substack{t\in \mathbb{P}:\\ t\leq (\ln x)^{1/14}}}\frac{\ln t}{t\varphi (t)}+ 2x  \sum_{\substack{t\in \mathbb{P}:\\ t\leq (\ln x)^{1/14}}} \frac{\ln t}{t\, \textup{exp} (b t^{-2} \sqrt{\ln (2x)})}\Biggr).\label{T8.BASIC.GENERAL}
\end{align} Since $\varphi (n) \geq c n/ \ln\ln(n+2)$ for any positive integer $n$, where $c$ is a positive absolute constant, we obtain
\[
\sum_{\substack{t\in \mathbb{P}:\\ t\leq (\ln x)^{1/14}}}\frac{\ln t}{t\varphi (t)} \leq \frac{1}{c}
\sum_{\substack{t\in \mathbb{P}:\\ t\leq (\ln x)^{1/14}}}\frac{\ln t \ln\ln (t+2)}{t^{2}}\leq
\frac{1}{c}
\sum_{n=1}^{\infty}\frac{\ln n \ln\ln (n+2)}{n^{2}}=c_1,
\]where $c_1$ is a positive absolute constant. We obtain
\begin{equation}\label{T8.BASIC.1}
\pi(2x)\sum_{\substack{t\in \mathbb{P}:\\ t\leq (\ln x)^{1/14}}}\frac{\ln t}{t\varphi (t)}\leq
c_1 \pi(2x)\leq c_1 a\,\frac{2x}{\ln (2x)}\leq 2c_1 a\,\frac{x}{\ln x} = c_2\,\frac{x}{\ln x},
\end{equation}where $c_2$ is a positive absolute constant.

For any prime $t \leq (\ln x)^{1/14}$, we have
\[
\frac{b \sqrt{\ln (2x)}}{t^{2}}\geq \frac{b \sqrt{\ln (2x)}}{(\ln x)^{1/7}}\geq
\frac{b (\ln x)^{1/2}}{(\ln x)^{1/7}}= b (\ln x)^{5/14}.
\] We obtain
\[
2x  \sum_{\substack{t\in \mathbb{P}:\\ t\leq (\ln x)^{1/14}}} \frac{\ln t}{t\, \textup{exp} (b t^{-2} \sqrt{\ln (2x)})} \leq 2x\,\textup{exp}\bigl(-b (\ln x)^{5/14}\bigr) \sum_{\substack{t\in \mathbb{P}:\\ t\leq (\ln x)^{1/14}}} \frac{\ln t}{t}.
\]

 Putting $k= \bigl[(\ln x)^{1/14}\bigr]$ and applying Lemma \ref{L:primes}, we have
 \[
 \sum_{\substack{t\in \mathbb{P}:\\ t\leq (\ln x)^{1/14}}} \frac{\ln t}{t}=
 \sum_{p\leq k} \frac{\ln p}{p}\leq c \ln k \leq c \ln \bigl((\ln x)^{1/14}\bigr)= \frac{c}{14}\ln\ln x=c_3\ln\ln x,
\]where $c_3 = c/14$ is a positive absolute constant.

We obtain
\[
2x  \sum_{\substack{t\in \mathbb{P}:\\ t\leq (\ln x)^{1/14}}} \frac{\ln t}{t\, \textup{exp} (b t^{-2} \sqrt{\ln (2x)})}\leq 2c_{3}x\,\textup{exp}\bigl(-b (\ln x)^{5/14}\bigr)\ln\ln x.
\]

Let us show that
\begin{equation}\label{T8.1.7.INEQ}
2c_{3}x\,\textup{exp}\bigl(-b (\ln x)^{5/14}\bigr)\ln\ln x \leq \frac{x}{\ln x}
\end{equation} or, that is equivalent,
\[
 2c_{3}\ln x \ln\ln x \leq \textup{exp}\bigl(b (\ln x)^{5/14}\bigr).
\]Taking logarithms, we obtain
\[
\ln(2c_3)+\ln\ln x + \ln\ln\ln x \leq b (\ln x)^{5/14}.
\] This inequality holds, if $x_0$ is chosen large enough. The inequality \eqref{T8.1.7.INEQ} is proved.

 We have
\begin{equation}\label{T8.BASIC.2}
2x  \sum_{\substack{t\in \mathbb{P}:\\ t\leq (\ln x)^{1/14}}} \frac{\ln t}{t\, \textup{exp} (b t^{-2} \sqrt{\ln (2x)})} \leq \frac{x}{\ln x}.
\end{equation} Substituting \eqref{T8.BASIC.1} and \eqref{T8.BASIC.2} into \eqref{T8.BASIC.GENERAL}, we obtain
\[
\sum_{\substack{t\in \mathbb{P}:\\ t\leq (\ln x)^{1/14}}} \frac{\lambda (q,t) \ln t}{t}\leq
C(E)(c_2 +1)\frac{x}{\ln x} = C_{1}(E) \frac{x}{\ln x},
\]where $C_{1}(E) >0$ is a constant depending only on $E$.

Applying \eqref{T8.N.A.ESTIM}, we obtain
\begin{align*}
\sum_{\substack{q\in \mathbb{P}:\\ \#E(\mathbb{F}_{q})< x}}\sum_{\substack{t\in \mathbb{P}:\\ t\leq (\ln x)^{1/14}}} \frac{\lambda (q,t) \ln t}{t}&\leq C_{1}(E)\,\frac{x}{\ln x}\,N_{A}(x)\\
 &\leq\frac{C_{1}(E)}{a_1}\,(N_{A}(x))^{2}= \gamma_{2}(E) (N_{A}(x))^{2},
\end{align*} where $\gamma_{2}(E) >0$ is a constant depending only on $E$. We see that the inequality \eqref{T6.INEQ3} with $\alpha=1/14$ holds.

By Theorem \ref{T6}, there are positive constants $c_{1}=c_{1}(\gamma_1)$ and $c_{2}=c_{2} (\gamma_{1}, \gamma_{2}(E), \alpha)$ such that
\[
\#\biggl\{1\leq n \leq x: r(n)\geq c_{1}\frac{N_{A}(x)}{\ln x}\biggr\}\geq c_{2}x\,\frac{N_{A}(x)}{N_{A}(x)+\rho_{A}(x)\ln x},
\]where
\[
r(n)=\#\{(p, q)\in \mathbb{P}^{2}: p+\#E(\mathbb{F}_q)=n\}.
\] Since $\alpha = 1/14$, $\gamma_1$ is a positive absolute constant, and $\gamma_2 (E)$ is a positive constant depending only on $E$, we see that $c_1$ is a positive absolute constant and $c_{2}$ is a positive constant depending only on $E$.

By \eqref{T8.N.A.ESTIM}, we have
\[
\#\biggl\{1\leq n \leq x: r(n)\geq c_{1}\frac{N_{A}(x)}{\ln x}\biggr\}\leq
\#\biggl\{1\leq n \leq x: r(n)\geq c_{1} a_{1}\frac{x}{(\ln x)^{2}}\biggr\}.
\] By \eqref{T8.rho.EST}, we have
\[
\rho_{A}(x)\ln x \leq 13 \sqrt{x}\ln x\leq a_{1}\frac{x}{\ln x},
\] if $x_0$ is chosen large enough (here $a_1$ is a positive absolute constant from \eqref{T8.N.A.ESTIM}). We obtain
\[
1 \leq \rho_{A}(x)\ln x \leq N_{A}(x),
\]and hence
\[
\frac{N_{A}(x)}{N_{A}(x)+\rho_{A}(x)\ln x}\geq
\frac{N_{A}(x)}{N_{A}(x)+ N_{A}(x)}= \frac{1}{2}.
\] We obtain
\[
\#\biggl\{1\leq n \leq x: r(n)\geq c_{1} a_{1}\frac{x}{(\ln x)^{2}}\biggr\}\geq
\frac{c_{2}(E)}{2}x.
\] We see that $c_{1} a_{1}$ is a positive absolute constant and $c_{2}(E)/2$ is a positive constant depending only on $E$. Let us denote $c_{1} a_{1}$ by $c_1$ and $c_{2}(E)/2$ by $c_{2}(E)$. Theorem \ref{T8} is proved.\\[-5pt]

\textsc{Proof of Theorem \ref{T9}.} We put
\[
A=\bigl\{a^{j^{b}}: j=0, 1, 2,\ldots\bigr\}.
\]We assume that $x\geq x_{0}(a, b)$, where $x_{0}(a, b) >0$ is a constant depending only on $a$ and $b$; we will choose the constant $x_{0}(a, b)$ later, and it will be large enough. We are going to apply Theorem \ref{T6}. It is clear that
\begin{gather*}
\textup{ord}_{A}(n)=\#\{j\geq 0: a^{j^{b}}=n\}\leq 1,\qquad n\in\mathbb{N},\\
N_{A}(x)=\#\bigl\{j\geq 0: a^{j^{b}}\leq x\bigr\} =y+1,\ \text{where $y:=\bigg[\biggl(\frac{\ln x}{\ln a}\biggr)^{1/b}\bigg]$}.
\end{gather*} Since $\textup{ord}_{A}(1)=1$, we see that
\[
\rho_{A}(x)= \max_{n\leq x} \textup{ord}_{A}(n) =1.
\] We have
\begin{equation}\label{T9.Proof.Ineq1}
y\geq 10,\qquad \biggl(\frac{\ln x}{\ln a}\biggr)^{1/b}\leq N_{A}(x)\leq 2 \biggl(\frac{\ln x}{\ln a}\biggr)^{1/b},\qquad N_{A}(x/2)\geq \frac{1}{\sqrt{8}}N_{A}(x),
\end{equation}if $x_{0}(a, b)$ is chosen large enough. It follows from \eqref{T9.Proof.Ineq1} that \eqref{T6.INEQ1} and \eqref{T6.INEQ2} hold.

 Let us show that
\begin{equation}\label{T9.SUPER.INEQ}
\sum_{0\leq k\leq y}\sum_{p\leq (\ln x)^{1/(2b)}}\frac{\lambda(k,p) \ln p}{p}\leq \gamma_{2}(a, b) y^{2},
\end{equation}where
\[
\lambda(k,p) = \#\Lambda (k, p),\qquad \Lambda(k,p):= \big\{k<j\leq y: a^{j^{b}}\equiv a^{k^{b}}\ \text{(mod $p$)}\},
\]and $\gamma_{2}(a, b)$ is a positive constant, depending only on $a$ and $b$.

 Fix an integer $k$ with $0\leq k \leq y$. We have
\[
\sum_{p\leq (\ln x)^{1/(2b)}} \frac{\lambda(k,p)\ln p}{p}=
\sum_{\substack{p\leq (\ln x)^{1/(2b)}:\\ p|a\text{ or }p|(a-1)}} +
\sum_{\substack{p\leq (\ln x)^{1/(2b)}:\\ (p,a)=1\text{ and }(p,a-1)=1}}  =S_1+ S_2.
\]Since $\lambda (k, p)\leq y$, we obtain
\[
S_1\leq y \sum_{\substack{p\leq (\ln x)^{1/(2b)}:\\ p|a\text{ or }p|(a-1)}}\frac{\ln p}{p}\leq y \sum_{\substack{p:\\ p|a\text{ or }p|(a-1)}}\frac{\ln p}{p}= c_{1}(a)y,
\]where $c_{1}(a)>0$ is a constant, depending only on $a$.

Let $p$ be in the range of summation of $S_2$. Since $(a,p)=1$, we have $a^{p-1}\equiv 1$ (mod $p$) (Fermat's theorem). Let $h_{a}(p)$ denote the order of $a$ modulo $p$, which is to say that $h_{a}(p)$ is the least positive integer $h$ such that $a^{h}\equiv 1\ \text{(mod $p$)}$. Since $(p,a-1)=1$, we see that $1<h_{a}(p)\leq p-1$. Let $j\in \Lambda (k, p)$. Then
\[
a^{j^{b}}\equiv a^{k^{b}}\ \text{(mod $p$)}.
\]Since
\[
(a^{k^{b}}, p)=1,
\]we obtain
\[
a^{j^{b}-k^{b}}\equiv 1\ \text{(mod $p$)}.
\]Hence (see, for example, \cite[Chapter VI]{Hardy.Wright}), $j^{b}-k^{b}\equiv 0$ (mod $h_{a}(p)$).

We use the following result of Konyagin \cite[Theorem 2]{Konyagin}. Let $m, n\in \mathbb{N}$. Let
\[
f(x)= \sum_{i=0}^{n} a_{i} x^{i},
\]where $a_{0}, \ldots, a_{n}$ are integers with $(a_{0},\ldots, a_{n}, m)=1$. By $\rho (f, m)$ we denote the number of solutions of the congruence $f(x)\equiv 0$ (mod $m$). Then
\begin{equation}\label{T9.K.EST}
\rho (f, m) \leq c n m^{1-1/n},
\end{equation}where $c$ is a positive absolute constant.

We define
\[
U=\big\{j\in \{0,\ldots, h_{a}(p)-1\}: j^{b}-k^{b}\equiv 0 \text{ (mod $h_{a}(p)$)}\big\}.
\] From \eqref{T9.K.EST} we obtain
\[
\#U\leq cb (h_{a}(p))^{1-1/b}.
\] It is easy to see that
\[
\lambda (k, p) \leq \sum_{j_{0}\in U} \#\bigl\{t\in \mathbb{Z}: 1\leq j_{0}+h_{a}(p)t\leq y\bigr\}.
\]We have
\[
(\ln x)^{1/ (2b)} \leq \biggl(\frac{\ln x}{\ln a}\biggr)^{1/b}-1,
\]if $x_{0}(a,b)$ is chosen large enough. Hence,
\[
h_{a}(p) \leq p-1 < p\leq (\ln x)^{1/ (2b)} \leq \biggl(\frac{\ln x}{\ln a}\biggr)^{1/b}-1\leq y.
\]Given $j_{0}\in U$, we have
\[
\#\{t\in \mathbb{Z}: 1\leq j_{0}+h_{a}(p)t\leq y\}\leq \frac{y}{h_{a}(p)}+1\leq \frac{2y}{h_{a}(p)}.
\] We obtain
\[
\lambda (k, p)\leq \frac{2y}{h_{a}(p)} \#U \leq \frac{2y}{h_{a}(p)} cb (h_{a}(p))^{1-1/b}=
c_{2}(b) \frac{y}{(h_{a}(p))^{1/b}},
\]where $c_{2}(b)=2cb$ is a positive constant, depending only on $b$.

We obtain
\begin{align}\label{T9.S2.LOW}
S_2 &= \sum_{\substack{p\leq (\ln x)^{1/(2b)}:\\ (p,a)=1\text{ and }(p,a-1)=1}}\frac{\lambda(k,p)\ln p}{p}\leq c_{2}(b) y\sum_{\substack{p\leq (\ln x)^{1/(2b)}:\\ (p,a)=1\text{ and }(p,a-1)=1}}\frac{\ln p}{p(h_{a}(p))^{1/b}}\notag\\
&\leq c_{2}(b) y\sum_{\substack{p:\\ (p,a)=1}}\frac{\ln p}{p(h_{a}(p))^{1/b}}.
\end{align}

We put
\[
D(z)=\sum_{n\leq z} d_n,\quad \text{where }d_n= \sum_{\substack{p:\\ (p,a)=1\\ h_{a}(p)=n}}\frac{\ln p}{p},\quad n=1, 2, \ldots.
\] Let $z\geq 100$. We have
\[
D(z)=\sum_{n\leq z} d_n= \sum_{n\leq z} \sum_{\substack{p:\\ (p,a)=1\\ h_{a}(p)=n}}\frac{\ln p}{p}.
\]Let $n$ and $p$ be in the range of summation. Then $h_{a}(p)=n$ and, hence, $a^{n}\equiv 1$ (mod $p$), i.e., $p| (a^{n}-1)$. We put $P(z)=\prod_{n\leq z}(a^{n}-1)$. We obtain
\[
D(z)\leq \sum_{p|P(z)}\frac{\ln p}{p} \leq c_{0}\ln\ln P(z),
\]where $c_{0}>0$ is an absolute constant. Since
\[
P(z)\leq \prod_{n\leq z} a^{n}= a^{1+2+\ldots+[z]}\leq a^{z^{2}},
\]we obtain $D(z)\leq c_1(a)\ln z$, where $c_1(a)>0$ is constant, depending only on $a$. Hence, $D(z)\leq c(a)\ln (z+1)$ for any real number $z\geq 1$, where $c(a)>0$ is a constant, depending only on $a$. Hence, $D(z)z^{-1/b}\to 0$ as $z\to +\infty$. Applying partial summation, we have
\[
\sum_{n\leq z}\frac{d_n}{n^{1/b}}=\frac{D(z)}{z^{1/b}}+ \frac{1}{b}\int_{1}^{z}\frac{D(t)}{t^{1+1/b}}\,dt
\] for any real $z\geq 1$. We obtain
\[
\sum_{n\geq 1}\frac{d_n}{n^{1/b}}= \frac{1}{b}\int_{1}^{+\infty}\frac{D(t)}{t^{1+1/b}}\,dt\leq
\frac{c(a)}{b}\int_{1}^{+\infty}\frac{\ln (t+1)}{t^{1+1/b}}\,dt = c_{3}(a,b),
\]where $c_3 (a, b)>0$ is a constant, depending only on $a$ and $b$.

We have
\begin{align*}
\sum_{n\geq 1}\frac{d_{n}}{n^{1/b}}&=\sum_{n\geq 1} \sum_{\substack{p:\\ (p,a)=1\\ h_{a}(p)=n}}\frac{\ln p}{p (h_{a}(p))^{1/b}}=\sum_{\substack{p:\\ (p,a)=1}} \frac{\ln p}{p (h_{a}(p))^{1/b}}\sum_{\substack{n\geq 1:\\ h_{a}(p)=n}}1\\
&=\sum_{\substack{p:\\ (p,a)=1}} \frac{\ln p}{p (h_{a}(p))^{1/b}}.
\end{align*} Hence (see \eqref{T9.S2.LOW}), $S_{2}\leq c_4(a,b)y$, where $c_4(a,b)>0$ is a constant, depending only on $a$ and $b$, and
\[
\sum_{p\leq (\ln x)^{1/(2b)}} \frac{\lambda(k,p)\ln p}{p}=S_1+S_2 \leq \bigl(c_1(a)+ c_4 (a, b)\bigr)y= c_{5}(a, b)y.
\]We obtain
\begin{align*}
\sum_{0\leq k\leq y}\sum_{p\leq (\ln x)^{1/(2b)}}\frac{\lambda(k,p)\ln p}{p}&\leq \sum_{0\leq k\leq y} c_{5}(a,b)y= c_{5}(a,b)y(y+1)\\
&\leq 2c_{5}(a, b)y^{2} = \gamma_{2}(a,b) y^{2},
\end{align*}where $\gamma_{2}(a,b)>0$ is a constant, depending only on $a$ and $b$. The inequality \eqref{T9.SUPER.INEQ} is proved. We see that \eqref{T6.INEQ3} holds with $\alpha=1/(2b)$.

 Applying Theorem \ref{T6}, we obtain
\begin{align*}
\#\{1&\leq n \leq x: \text{there are $p\in \mathbb{P}$ and $j\in \mathbb{Z}_{\geq 0}$ such that $p+ a^{j^{b}}=n$}\}\\
 &\geq c_{2}(a,b)x\,\frac{N_{A}(x)}{N_{A}(x)+\rho_{A}(x)\ln x},
\end{align*}where $c_{2}(a,b)>0$ is a constant, depending only on $a$ and $b$.

We have
\[
\biggl(\frac{\ln x}{\ln a}\biggr)^{1/b}\leq \ln x,
\] if $x_{0}(a,b)$ is chosen large enough. Applying \eqref{T9.Proof.Ineq1} and taking into account that $\rho_{A}(x)=1$, we obtain
\begin{align*}
\frac{N_{A}(x)}{N_{A}(x)+\rho_{A}(x)\ln x}&\geq \frac{(\ln x/\ln a)^{1/b}}{2(\ln x/\ln a)^{1/b}+\ln x}\\
&\geq
\frac{(\ln x/\ln a)^{1/b}}{2\ln x+\ln x}=
\frac{1}{3 (\ln a)^{1/b}}\frac{1}{(\ln x)^{1-1/b}}.
\end{align*} We obtain
\begin{align*}
\#\{1&\leq n \leq x: \text{there are $p\in \mathbb{P}$ and $j\in \mathbb{Z}_{\geq 0}$ such that $p+ a^{j^{b}}=n$}\}\\
 &\geq \frac{c_{2}(a,b)}{3 (\ln a)^{1/b}} \frac{x}{(\ln x)^{1-1/b}} = c_{1}(a, b)\frac{x}{(\ln x)^{1-1/b}}  ,
\end{align*} where $c_{1}(a,b)>0$ is a constant, depending only on $a$ and $b$.

We have
\begin{align*}
\#\{1&\leq n \leq x: \text{there are $p\in \mathbb{P}$ and $j\in \mathbb{Z}_{\geq 0}$ such that $p+ a^{j^{b}}=n$}\}\\
&\leq \#\{(p,j)\in \mathbb{P}\times \mathbb{Z}_{\geq 0}: p\leq x\text{ and }a^{j^{b}}\leq x\}= \#\{p: p\leq x\} \#\{j\geq 0: a^{j^{b}}\leq x\}\\
&=\pi(x) N_{A}(x)\leq c\,\frac{x}{\ln x}\,2\,\biggl(\frac{\ln x}{\ln a}\biggr)^{1/b}= c_{2}(a, b)\, \frac{x}{(\ln x)^{1-1/b}},
\end{align*}where $c_{2}(a, b)>0$ is a constant, depending only on $a$ and $b$.

Since
\[
3=p+a^{j^{b}}
\]for $p=2$ and $j=0$, the claim in the case $3 \leq x < x_{0}(a,b)$ is trivial. Theorem \ref{T9} is proved.

\section{Acknowledgements}

I would like to thank Sergei Konyagin for many useful conversations and suggestions. Also I would like to thank the anonymous referee for many useful comments.


\begin{thebibliography}{99}

  \bibitem{David.Wu}
   C.~David and J.~Wu, Pseudoprime reductions of elliptic curves, \emph{Canad. J. Math.}, \textbf{64} (2012), no. 1,  81--101.


   \bibitem{Hardy.Wright}
  G.\,H.~Hardy and E.\,M.~Wright, \emph{An introduction to the theory of numbers}, 6th ed., Oxford Univ. Press, Oxford, 2008, xxii+621 pp.

   \bibitem{Hasse}
  H.~Hasse, Abstrakte Begr\"{u}ndung der komplexen Multiplikation und Riemannsche Vermutung in Funktionenk\"{o}rpern, \emph{Abh. Math. Sem. Hamburg}, \textbf{10} (1934), 325--348.

   \bibitem{Konyagin}
 S.\,V.~Konyagin, On the number of solutions of an $n$th degree congruence with one unknown, \emph{Sb. Math.}, \textbf{37} (1980), no. 2, 151--166.

   \bibitem{Maynard}
  J.~Maynard, Dense clusters of primes in subsets, \emph{Compos. Math.}, \textbf{152} (2016), no. 7, 1517--1554.


   \bibitem{Prachar}
  K.~Prachar, \emph{Primzahlverteilung}, Springer-Verlag, Berlin--G\"{o}ttingen--Heidelberg, 1957, x+415 pp.


    \bibitem{Murty}
   M.~Ram Murty, \emph{Problems in analytic number theory}, 2nd ed., Grad. Texts in Math., \textbf{206}, Readings in Math., Springer, New York, 2008, xxii+502 pp.


 \bibitem{Romanoff}
  N.\,P.~Romanoff, \"{U}ber einige S\"{a}tze der additiven Zahlentheorie, \emph{Math. Ann.}, \textbf{109} (1934), 668--678.


 \bibitem{Schnirelmann}
  L.~Schnirelmann, \"{U}ber additive Eigenschaften von Zahlen, \emph{Math. Ann.}, \textbf{107} (1933), 649--690.






  \end{thebibliography}
\end{document}